\documentclass[a4paper, 11pt]{article} 
\usepackage{amsfonts,amsmath,amssymb,amscd,color,stmaryrd,latexsym}
\usepackage[all]{xy} 
\usepackage[T1]{fontenc}

\def\End{\mathop{\rm End}\nolimits}

\def\vol{\mathop{vol}\nolimits}

\def\Re{\mathop{\rm Re}\nolimits}

\def\deg{\mathop{\rm deg}\nolimits}

\def\cliff{\mbox{\sl Cliff}}

\def\ric{\mathop{\rm Ric}\nolimits}
\def\mod{\mathop{\rm mod}\nolimits}
\def\note#1{\marginpar{\raggedright\if@twoside\ifodd\c@page\raggedleft\fi\fi\sf\scriptsize RMK: #1}}

\newcommand{\mf}{\mathfrak}
\newcommand{\mb}{\mathbf}
\newcommand{\mc}{\mathcal}

\newcommand{\R}{\mathbb{R}}
\newcommand{\Z}{\mathbb{Z}}
\newcommand{\C}{\mathbb{C}}

\newcommand{\eunorm}[1]{\parallel\mbox{\hspace{-3pt}} #1\mbox{\hspace{-3pt}}\parallel}

\renewcommand{\iff}{if and only if }
\newcommand{\wrt}{with respect to }

\newcommand{\Id}{\mathrm{Id}}

\newcommand{\gtilde}{\widetilde{\mc{G}}}
\newcommand{\htimes}{\widehat{\otimes}}
\newcommand{\D}{\mathrm{D}}

\newenvironment{rmk}{\begin{trivlist}\item[]{\bf Remark:}\setlength{\parindent}{0pt}}{\end{trivlist}}
\newenvironment{ex}{\begin{trivlist}\item[]{\bf Example:}\setlength{\parindent}{0pt}}{\end{trivlist}}
\newenvironment{prf}{\begin{trivlist}\item[]{\bf Proof: }}{\hfill $\blacksquare$ \end{trivlist}}

\newtheorem{thm}{Theorem}[section]
\newtheorem{definition}[thm]{Definition}
\newtheorem{prp}[thm]{Proposition}
\newtheorem{lem}[thm]{Lemma}
\newtheorem{cor}[thm]{Corollary}

\begin{document}
\thispagestyle{empty}
\begin{flushright}
\hfill{MPP-2005-107}
\end{flushright}
\vspace{12pt}

\begin{center}{\LARGE{\bf Generalised geometries, constrained critical points and Ramond--Ramond fields}}

\vspace{35pt}

{\bf Claus Jeschek}$^a$ and {\bf Frederik Witt}$^b$

\vspace{15pt}

$^a$ {\it  Max--Planck--Institut f\"ur Physik,\\
D--80805 M\"unchen, F.R.G.\\
e--mail: {\sf jeschek@mppmu.mpg.de}}\\[1mm]

\vspace{15pt}

$^b$ {\it Freie Universit\"at Berlin, FB Mathematik und Informatik,\\
Arnimallee 3, D-14195 Berlin, F.R.G.\\e--mail: {\sf fwitt@math.fu-berlin.de}}\\[1mm]

\vspace{40pt}

{\bf ABSTRACT}
\end{center}
In the context of generalised geometry we investigate reductions to $SU(m)\times SU(m)$ together with an integrability condition which in dimension $6$ describes the geometry of type II supergravity compactifications.

\bigskip


\smallskip

\section{Introduction}
Generalised geometries were introduced by Hitchin in~\cite{hi03} and have played an important r\^ole in string theory ever since. In a nutshell, we associate with a closed $3$--form $H$ on $M^n$ (the $H$--{\em flux} in physicists' language) a vector bundle $\mb{E}\to M$ with structure group $SO(n,n)$. 
The space of sections carries a natural bracket -- the Courant bracket. More importantly for us, $\mb{E}$ admits a natural spin structure, and $H$ induces an elliptic complex on the spinor fields. These can be thought of as even or odd forms with differential defined by $d_H\alpha=d\alpha+H\wedge\alpha$. In analogy with classical geometries described by reductions inside the $GL(n)$--frame bundle and integrability conditions involving the Lie bracket or the exterior differential, we can consider reductions inside $SO(n,n)$ and integrability conditions involving the Courant bracket or $d_H$. In particular, the usual classical geometries (such as complex, K\"ahler or $G_2$--manifolds) have natural counterparts in this setting (cf.~\cite{hi03},~\cite{gu04},~\cite{wi05}). 

\bigskip

In this article we shall be concerned with {\em generalised} $SU(m)$--structures~\cite{gmwi06}, i.e.\ 
reductions from $\R_{\not=0}\times Spin(2m,2m)$ to $SU(m)\times SU(m)$. While a usual $SU(m)$--structure (i.e.\ an almost Calabi--Yau) can be characterised by a(n almost) symplectic $2$--form and a complex $m$--form subject to various compatibility conditions, generalised $SU(m)$--structures are characterised by mixed degree complex forms $\rho_0$ and $\rho_1$ of special type, together with a scaling function, the so--called {\em dilaton}. More precisely, a generalised $SU(m)$--structure gives rise to a spinnable metric $g$ with spinor bundle $\Delta$ and two chiral spinor fields $\Psi_{L,R}\in\Gamma(\Delta)$. Under the isomorphism $\Omega^*(M)\otimes\C\cong\Delta\otimes\Delta$, we get (up to a $B$--field transform)
\begin{equation}\nonumber
\rho_0=e^{-\phi}\mc{A}(\Psi_L)\otimes\Psi_R,\quad\rho_1=e^{-\phi}\Psi_L\otimes\Psi_R,
\end{equation}
where $\mc{A}$ is the conjugate--linear, $Spin(2m)$--equivariant {\em charge conjugation operator}. As we will show, both $\rho_0$ and $\rho_1$ give rise to a(n almost) generalised Calabi--Yau structure in the sense of Hitchin~\cite{hi03}. Further, our notion of a generalised $SU(m)$--structure coincides with Gualtieri's notion of a(n almost) generalised Calabi--Yau metric~\cite{gu04}.

\bigskip

Since we have a natural metric, we can perform Hodge theory for the above mentioned elliptic complex. In particular, the real parts $\Re(\rho_k)$ are harmonic \iff  
\begin{equation}\label{harmonicint}
d_H\rho_0=0,\quad d_H\rho_1=0.
\end{equation} 
In this case, $\rho_{0,1}$ define integrable generalised Calabi--Yau structures, while the pair $(\rho_0,\rho_1)$ gives an integrable generalised Calabi--Yau metric. Moreover we show in~\cite{jewi05} that a generalised $SU(m)$--structure with~(\ref{harmonicint}) is equivalent to the geometry of type II compactifications to $6$ dimensions with vanishing Ramond--Ramond fields. As a consequence, the underlying spinor fields $\Psi_L$ and $\Psi_R$ are parallel with respect to the spin connections $\nabla^{\pm}_X=\nabla^g_X\pm(X\llcorner H/4)\cdot$, where $\nabla^g$ denotes the Levi--Civita connection. This type of spinor field equations, however, implies $H=0$ for $M$ {\em compact}, so that the holonomy of the metric must be contained in $SU(m)$. This stands in striking contrast to {\em generalised K\"ahler manifolds} \big(whose structure group is $U(m)\times U(m)$\big) for which genuinely generalised examples do exist~\cite{cagu07},~\cite{hi06}. 

\bigskip

We therefore consider the general ansatz
\begin{equation}\label{genintcond}
d_H\rho_0=F_0,\quad d_H\rho_1=F_1
\end{equation}
for (complex valued) forms of mixed degree $F_{0,1}$. Equations of this type also occur in generalised calibration theory, where calibration forms satisfying~(\ref{genintcond}) give rise to calibrated submanifolds minimising the so-called Wess--Zumino term of the D--brane energy~\cite{gmwi06}. If $F_{0,1}$ are of special algebraic type, then the underlying spinor fields satisfy
\begin{eqnarray*}
\nabla_X^+\Psi_L & = & e^{\phi}\big(a_L\overline{F_0}\cdot X\cdot\Psi_R+b_LF_1\cdot X\cdot\mc{A}(\Psi_R)\big)\nonumber\\
\nabla_X^-\Psi_R & = & e^{\phi}\big(a_R\widehat{F_0}\cdot X\cdot\Psi_L+b_R\widehat{F_1}\cdot X\cdot\mc{A}(\Psi_L)\big),
\end{eqnarray*}
where $\widehat{\cdot}\,$ is a sign--twisting operator on forms and $a_{L,R}$ and $b_{L,R}$ are constants depending on $m$. We refer to these $F_{0,1}$ of special type as {\em Ramond--Ramond fields}, since a generalised $SU(3)$--structure with~(\ref{genintcond}) for $F_{0,1}$ either real or purely imaginary accounts for the Ramond--Ramond fields in type II compactifications. Note that similar statements have been asserted in the physics literature by Gra\~na et al.~\cite{gmpt04},~\cite{gmpt05}, though their findings in~\cite{gmpt05} seem to contradict the special case considered in~\cite{gmpt04}. Simple, but compact examples are provided by ordinary $SU(3)$--structures of type $W_2^{\pm}$ or $W_3$, cf.~\cite{chsa02}.  From a mathematical point of view, the condition~(\ref{genintcond}) can be interpreted as the Euler--Lagrange equation of the a constrained variational problem if the Ramond--Ramond fields are $d_H$--exact (as required by physics).

\bigskip

Finally, a comment on conventions. In this article, we refer to any reduction of the structure group to a subgroup $G$ as a $G$--{\em structure}. This is purely topological and does not imply any sort of integrability. For instance, in our jargon a complex structure on $M^{2m}$ is simply an almost complex manifold in the sense that the structure group reduces to $GL(m,\C)$. Though we sometimes use the additional qualifier ``almost'' we usually omit it for sake of simplicity. Bearing this caveat in mind, no confusion should appear.

\bigskip

The paper is organised as follows. In Section~\ref{gengeo} we briefly recall the setup of generalised geometry following~\cite{hi03},~\cite{hi05} and~\cite{gmwi06}. In Section~\ref{gencom}, we explain the reduction mechanism which leads to generalised $SU(m)$--structures and consider the relation with generalised Calabi--Yau structures and metrics. We interlude in Section~\ref{rrfields} and discuss supergravity compactifications to $6$ dimensions and the r\^ole played by the Ramond--Ramond fields. We generalise the resulting integrability condition to any generalised $SU(m)$--structure in Section~\ref{varprinc}. We analyse this condition from a Hodge theoretic point of view and translate it to the setting of generalised geometry. 
%
%
%
%
%
\section{Generalised geometries}\label{gengeo}
%
%
%
\subsection{{\bf The group $\mb{SO(n,n)}$}}\label{groupnn}
We briefly recall the algebraic setup of generalised geometry as introduced by Hitchin~\cite{hi03}.

\bigskip

Let $(\R^{n,n},g)$ or simply $\R^{n,n}$ denote the vector representation of $O(n,n)$ with invariant inner product $g$. Because of split signature, we can always choose a decomposition $\R^{n,n}=W\oplus W'$ into two maximally isotropic subspaces $W$ and $W'$, i.e.\ $g_{|W,W'}\equiv0$ and $\dim W,W'=n$. Identifying the isotropic complement $W'$ with $W^*$, we obtain an isometry $(\R^{n,n},g)\cong\big(W\oplus W^*,(\cdot\,,\cdot)\big)$
with $(w,\xi)=\xi(w)/2$ for $w\in W$ and $\xi\in W^*$. Moreover, this induces a canonic orientation on $\R^{n,n}$, so that we can think of $W\oplus W^*$ as an $SO(n,n)$--representation space. Note that the choice of $W$ gives rise to the subgroup $GL(W)\cong GL(n)\subset O(n,n)$, as for $A\in GL(d)$, $(Aw,A^*\xi)=(w,\xi)$. The two connected components of $GL(n)$ single out the two connected components of $O(n,n)$ which make up the group $SO(n,n)$. Furthermore, an orientation on $W$ gives rise to a space-- and timelike orientation which is preserved by the identity component $SO(n,n)_+$. 

\bigskip

The Lie algebra of $SO(n,n)_+$ is, regarded as being acted on by the subgroup $GL(n)_+$, the direct sum of the irreducible pieces
\begin{equation}\label{liealsplit}
\mf{so}(n,n)\cong\Lambda^2(W\oplus W^*)=\Lambda^2W\oplus W\otimes W^*\oplus\Lambda^2W^*.
\end{equation}
The subspace $W\otimes W^*=\mf{gl}(n)$ is just the Lie algebra of $GL(n)_+$ which shows that the inclusion of $GL(n)_+$ into $SO(n,n)_+$ lifts to an inclusion of $GL(n)_+$ into $Spin(n,n)_+$. The spin representations $S_{\pm}$ of $Spin(n,n)$ can be canonically constructed from the choice of $W$. We define for $X\oplus\xi\in W\oplus W^*$ an action on $\rho\in\Lambda^*W^*$ by
\begin{equation}\nonumber
(X\oplus\xi)\bullet\rho=-X\llcorner\rho+\xi\wedge\rho.
\end{equation}
It is immediate from the definition that this squares to minus the identity on unit vectors, so $\cliff(W\oplus W)\cong\End(\Lambda^*W^*)\cong\cliff(n,n)$. Moreover, the action of $\cliff(n,n)^{ev}$ preserves the parity of elements in $\Lambda^*$, i.e.\ even or odd forms are again mapped to even or odd forms. Therefore, $S_{\pm}=\Lambda^{ev,od}W^*$ are the two irreducible spin representations of $Spin(n,n)$. As a $GL(n)_+$--module, we find $A\bullet\rho=\sqrt{\det A}\cdot A^*\rho$, where $A^*\rho$ is the standard action of $A\in GL(n)_+$ on $\rho\in\Lambda^{ev,od}W^*$. Hence, there is a $GL(n)_+$--equivariant isomorphism $S_{\pm}\stackrel{\cong}{\to}\Lambda^{ev,od}W^*\otimes\sqrt{\Lambda^nW}$. The choice of an $n$--vector $\nu$ allows us to safely identify spinors with forms. We shall write this isomorphism as $\rho\mapsto\rho^{\nu}$. Note that $\nu_1=e^{2\phi}\nu_0$ implies $\rho^{\nu_1}=e^{-\phi}\rho^{\nu_0}$. Next let $[\,\cdot\,]^n$ denote projection on the top degree component of $\alpha\in\Lambda^*W^*$, and $\widehat{\,\cdot\,}$ the sign--changing operator defined on forms of degree $p$ by 
\begin{equation}\nonumber
\widehat{\alpha^p}=(-1)^{p(p+1)/2}\alpha^p.
\end{equation}
A $Spin(n,n)_+$--invariant inner product on $\Gamma(S)$ is given by $\langle\rho,\tau\rangle=\nu\big([\rho^{\nu}\wedge\widehat{\tau^{\nu}}]^n\big)\in\R$ (clearly, this does not dependent on the choice of $\nu$).

\bigskip

A part from the $GL(n)_+$--action, a further class of transformations deserves special interest, namely the so--called $B$--{\em field transformations}. As we also see from~(\ref{liealsplit}), any 2--form $B=\sum B_{kl}w^k\wedge w^l\in\Lambda^2W^*$ can be regarded as an element in $\mf{so}(n,n)$. It therefore acts through exponentiation as an element of both $Spin(n,n)_+$ and $SO(n,n)_+$. Explicitly, we map $B$ into $\cliff(n,n)$ via $B\mapsto\sum_{kl}B_{kl}w^k\bullet w^l$, so it acts on spinors by
\begin{equation}\nonumber
e^B\bullet\rho=(1+B+\frac{1}{2}B\bullet B+\ldots)\bullet\rho=(1+B+\frac{1}{2}B\wedge B+\ldots)\wedge\rho=e^B\wedge\rho.
\end{equation}
On the other hand, $B$ becomes a skew--symmetric linear operator $W\to W^*$ as an element of $\Lambda^2(W\oplus W^*)$. Under
the identification $\zeta\wedge\xi(X)=(\zeta,X)\xi-(\xi,X)\zeta=X\llcorner(\zeta\wedge\xi)/2$, we obtain
\begin{equation}\nonumber
e^B_{SO(n,n)}(X\oplus\xi)=\left(\begin{array}{cc} 1\,\, & 0\\ B/2\,\, & 1\end{array}\right)\left(\begin{array}{c} X\\ \xi\end{array}\right)
\end{equation}
on $W\oplus W^*$, where $B(X)=X\llcorner B$.
%
%
%
\subsection{The generalised tangent bundle}\label{gentangbun}
Next, we consider the global situation (cf.~\cite{gmwi06} and~\cite{hi05} for details).

\bigskip

Consider a triple $(M,\nu,H)$ consisting of an orientable manifold $M$ of dimension $n$, a nowhere vanishing $n$--vector field $\nu$ \big(i.e.\ a section $\nu\in\Gamma(\Lambda^nTM)$\big) and a closed $3$--form $H$ (the $H$--{\em flux} in physicists' terminology). Choose a convex cover of coordinate neighbourhoods $\{U_a\}$ of $M$, whose induced transition functions of the tangent bundle are $s_{ab}:U_{ab}=U_a\cap U_b\to GL(n)_+
 $. Locally, $H_{|U_a}=dB^{(a)}$ for some $B^{(a)}\in\Omega^2(U_a)$. Then we can define the closed $2$--forms $\beta^{(ab)}=B^{(a)}_{|U_{ab}}-B^{(b)}_{|U_{ab}}\in\Omega^2(U_{ab})$. We twist the transition functions
\begin{equation}\nonumber
S_{ab}=\left(\begin{array}{cc} s_{ab} & 0\\ 0 & (s_{ab}^{-1})^{\rm tr}\end{array}\right)\in GL(n)_+
\subset SO(n,n)_+
\end{equation}
of $TM\oplus T^*\!M$ to get new transition functions $\sigma_{ab}=S_{ab}\circ e^{2\beta^{(ab)}}:U_{ab}\to SO(n,n)_+
$. These give rise to the {\em generalised tangent bundle}
\begin{equation}\nonumber
\mb{E}=\mb{E}(H)=\prod U_a\times(\R^n\oplus\R^{n*})/\sim_{\sigma_{ab}},
\end{equation}
which we can view as an extension of $T^*\!M$ by $TM$. As such it is isomorphic, as a $C^{\infty}$--vector bundle, to $TM\oplus T^*\!M$. In particular, we can identify any section $X\oplus\xi\in\Gamma(TM\oplus T^*\!M)$, locally given by $X_a\oplus\xi_a$, with the section of $\mb{E}$ locally given by $e^{2B^{(a)}_a}(X_a\oplus\xi_a)=X_a\oplus B^{(a)}_a(X_a)+\xi_a$.  

\bigskip

The importance of the $H$--field appears when we discuss the natural elliptic complex on spinor fields associated with $\mb{E}$. A canonic spin structure is specified by the transition functions $\widetilde{\sigma}_{ab}=\widetilde{S}_{ab}\bullet e^{\beta^{(ab)}_b}$, where $\widetilde{S}_{ab}\in GL(n)_+\subset Spin(n,n)_+$ denotes the lift of $S_{ab}$ and $\beta^{(ab)}$ is exponentiated to $Spin(n,n)_+$. The chiral spinor bundles associated with $\mb{E}$ are
\begin{equation}\nonumber
\mb{S}_{\pm}=\mb{S}(\mb{E})_{\pm}=\coprod\limits_aU_a\times S_{\pm}/\sim_{\widetilde{\sigma}_{ab}},\quad\mb{S}=\mb{S}(\mb{E})_{\pm}=\mb{S}_+\oplus\mb{S}_-.
\end{equation}
An $\mb{E}$--spinor field $\rho$ is thus represented by a collection of smooth maps $\rho_a:U_a\to S_{\pm}$ with $\rho_a=\widetilde{\sigma}_{ab}\bullet\rho_b$. Fix a trivialisation $\nu_0$ of $\Lambda^n\R^n$ so that $\nu_{|U_a}=\nu_a=\lambda_a^{-2}\nu_0\in\Lambda^n\R^n$ with $\lambda_a^{-2}\in C^{\infty}(M)$. The map $\rho_a\mapsto e^{-B^{(a)}_a}\wedge\lambda_a\cdot\rho_a$ glues to a well--defined isomorphism $\rho\in\Gamma(\mb{S}_{\pm})\to\rho^{\nu}\in\Omega^{ev,od}(M)$ with $\big((X\oplus\xi)\bullet\rho\big)^{\nu}=-X\llcorner\rho^{\nu}+\xi\wedge\rho^{\nu}$. In the same vein, the $Spin(n,n)_+$--invariant form $\langle\cdot\,,\cdot\rangle$ induces on $\Gamma(\mb{S})$ a globally defined inner product $\langle\rho,\tau\rangle=\nu([\rho^{\nu}\wedge\widehat{\tau^{\nu}}]^n)$.

\bigskip

The choice of $\nu$ gives also rise to the differential operator
\begin{equation}\nonumber
d_{\nu}:\Gamma\big(\mb{S}_{\pm}\big)\to\Gamma\big(\mb{S}_{\mp}\big),\quad(d_{\nu}\rho)_a=\lambda^{-1}_a\cdot d_a(\lambda_a\cdot\rho_a),
\end{equation}
where $d_a$ is the usual differential applied to forms $U_a\to\Lambda^*\R^{n*}$, i.e.\ $(d\alpha)_a=d_a\alpha_a$. Then
\begin{equation}\nonumber
(d_{\nu}\rho)^{\nu}=d_H\rho^{\nu}
\end{equation}
where $d_H$ acts on even or odd forms via
\begin{equation}\nonumber
d_H\alpha=d\alpha+H\wedge\alpha.
\end{equation} 
Since $d_H$ is just the ordinary exterior differential $d$ up to the $0$--order term $H\wedge$, $d_H$ (and thus $d_{\nu}$) defines an elliptic complex on $\Omega^*(M)=\Omega^{ev}(M)\oplus\Omega^{od}(M)$ which computes the so--called {\em twisted cohomology} of $M$. Up to isomorphism, the cohomology groups $\mb{H}^{ev,od}(M,H)$ only depend on $[H]$. We denote by $\mb{H}^{\pm}(\mb{E})$ the corresponding cohomology induced by $d_{\nu}$ (which up to natural isomorphism does not depend on $\nu$).
%
%
%
\subsection{Generalised Riemannian metrics}\label{genriemmet}
\begin{definition} 
A {\em generalised (Riemannian) metric} for the generalised tangent bundle $\mb{E}$ is the choice of a maximally positive definite subbundle $\mb{V}^+$.
\end{definition}

The choice of $\mb{V}^+$ is complemented by a maximally negative definite subbundle $\mb{V}^-$ so that in terms of structure groups, a generalised metric corresponds to the reduction of $SO(n,n)_+$ to $SO(\mb{V}^+)\times SO(\mb{V}^-)\cong SO(n,0)\times SO(0,n)$. A generalised metric is thus equivalent to a splitting of the exact sequence $0\to T^*\!M\to\mb{E}\to TM\to 0$. We denote by $X^{\pm}$ the lift of vector fields $X\in\Gamma(TM)$ to sections in $\Gamma(\mb{V}^{\pm})$. Locally, $X^{\pm}$ corresponds to smooth maps $X_a^{\pm}:U_a\to\R^n\oplus\R^{n*}$ with $X^{\pm}_a=\sigma_{ab}(X^{\pm}_b)$, and $X^{\pm}_a=X_a\oplus P^{\pm}_aX_a$ for linear isomorphisms $P^{\pm}_a:\R^n\to\R^{n*}$. The symmetric parts of the isomorphisms $P_a^+$ patch together to give a globally defined Riemannian metric (using the lifts $X^-$ yields the same metric). The skew--symmetric parts yield a globally defined $2$--form $B$. Passing to $\mb{E}(H+dB)$ if necessary, we may assume that $B=0$. Hence a generalised metric on $M$ boils down to a pair $(H,g)$. Further, we have the canonical Riemannian $n$--vector $\nu_g$. We shall simply write $\rho$ for $\rho^{\nu_g}$ if there is no danger of confusion.

\bigskip

A generalised metric induces also further structure on spinor fields~\cite{gmwi06}. Since $M$ is orientable, so is $\mb{V}^{\pm}$. Hence, we obtain a Riemannian volume element $\varpi_{\mb{V}^-}$ on each fibre. Let
\begin{equation}\nonumber
\widetilde{\alpha}=(-1)^p\alpha
\end{equation}
for $\alpha\in\Omega^p(M)$. Recast in terms of forms, the action of the operator $\gtilde=\varpi_{V^-}\bullet$ on $\mb{S}(\mb{E})$ is
\begin{equation}\nonumber
(\gtilde\rho)^{\nu}=\left\{\begin{array}{ll}$n$\mbox{ even:}\quad&\star_g\widehat{\rho^{\nu}}\\$n$\mbox{ odd:}\quad&\star_g\widehat{\widetilde{\rho^{\nu}}}\end{array}\right.=:\gtilde\rho^{\nu},
\end{equation}
where $\star_g$ denotes the Hodge operator associated with the metric. We have the identities 
\begin{equation}\nonumber
\gtilde^2=(-1)^{n(n+1)/2}Id\quad\mbox{ and}\quad\langle\gtilde\rho,\tau\rangle=(-1)^{n(n+1)/2}\langle\rho,\gtilde\tau\rangle.
\end{equation} 
In particular, $\widetilde{\mc{G}}$ is an isometry for $\langle\cdot\,,\cdot\rangle$ and defines a complex structure on $\mb{S}$ if $n\equiv1,2\mod4$. Finally, we define the bilinear form (with $[n/2]=m$ for $n=2m,\,2m+1$) 
\begin{equation}\nonumber
\mc{Q}^{\pm}(\rho,\tau)=\pm(-1)^m\langle\rho,\gtilde\tau\rangle=g(\rho,\tau)\mbox{ on }\mb{S}_{\pm},\quad m=[\frac{n}{2}].
\end{equation}

\bigskip

If in addition $M$ is spinnable, the presence of a generalised metric also implies a useful description of the complexified spinor modules $\mb{S}_{\pm}\otimes\C$. Let $\Delta_n=\Delta$ be the $Spin(n)=Spin(n,0)$--module of {\em (Dirac) spinors}. It is irreducible for $n$ odd, while for $n$ even $\Delta_n=\Delta_+\oplus\Delta_-$, where $\Delta_{\pm}$ are the irreducible representation spaces of {\em chiral} spinors. As a complex vector space, $\Delta_n\cong\C^{2^{[n/2]}}$. Further, $\Delta_n$ carries an hermitian inner product $q$ (which we take to be conjugate--linear in the first argument) for which $q(a\cdot\Psi,\Phi)=q(\Psi,\widehat{a}\cdot\Phi)$, $a\in\cliff(TM\otimes\C,g^{\C})$. There exists a conjugate--linear, $Spin(n)$--equivariant endomorphism $\mc{A}$ of $\Delta_n$ (the {\em charge conjugation operator} in physicists' language) such that~\cite{wa89}
\begin{equation}\label{Aequiv}
\mc{A}(X\cdot\Psi)=(-1)^{m+1}X\cdot \mc{A}(\Psi)\mbox{ and }\mc{A}^2=(-1)^{m(m+1)/2}Id,\quad m=[\frac{n}{2}].
\end{equation}
The operator $\mc{A}$ preserves chirality \iff $n=2m$ with $m$ even. The $Spin(n)$--invariant bilinear form $\mc{A}(\Psi,\Phi)=q(\mc{A}(\Psi),\Phi)$ satisfies
\begin{equation}\nonumber
\mc{A}(\Psi,\Phi)=(-1)^{m(m+1)/2}\mc{A}(\Phi,\Psi)\mbox{ and }\mc{A}(X\cdot\Psi,\Phi)=(-1)^m\mc{A}(\Psi,X\cdot\Phi),\quad m=[\frac{n}{2}].
\end{equation}
We can embed the tensor product $\Delta_n\otimes\Delta_n$ into $\Lambda^*\C^n$ by sending $\Psi_L\otimes\Psi_R$ to the form of mixed degree 
\begin{equation}\nonumber
[\Psi_L\otimes\Psi_R](X_1,\ldots,X_n)=\mc{A}\big(\Psi_L,(X_1\wedge\ldots\wedge X_n)\cdot\Psi_R\big).
\end{equation} 
In fact, this is an isomorphism for $n$ even. In the odd case, we obtain an isomorphism by concatenating $[\cdot\,,\cdot]$ with projection on even or odd forms, which we write as $[\cdot\,,\cdot]^{ev,od}$. Since this map is $Spin(n)$--equivariant, it acquires global meaning over $M$, and we use the same symbol for the resulting map $\Gamma(\Delta\otimes\Delta)\to\Omega^{ev,od}(M)$ (referred to as the {\em fierzing} map in the physics' literature). Let
\begin{equation}\nonumber
[\cdot\,,\cdot]^{\mc{G}\,ev,od}:\Gamma(\Delta\otimes\Delta)\stackrel{[\cdot,\cdot\,]^{ev,od}}{\longrightarrow}\Omega^{ev,od}(M)\otimes\C\stackrel{\nu_g}{\longrightarrow}\Gamma\big(\mb{S}_{\pm}\otimes\C\big).
\end{equation}
Clifford action on $\Delta$ and $\mb{S}_{\pm}$ relate via
\begin{equation}\label{commut}
\begin{array}{lcl}
\,[X\cdot\Psi_L\otimes\Psi_R]^{\mc{G}} & = & (-1)^{n(n-1)/2}X^+\bullet[\Psi_L\otimes\Psi_R]^{\mc{G}}\\[4pt]
\,[\Psi_L\otimes Y\cdot\Psi_R]^{\mc{G}} & = & -Y^-\bullet\widetilde{[\Psi_L\otimes\Psi_R]^{\mc{G}}}.
\end{array}
\end{equation}
As $a_1\cdot\ldots\cdot a_{2l}\in Spin(0,n)$ acts on $\Psi\in\Delta_n$ via $(-1)^la_1\cdot\ldots\cdot a_{2l}\cdot\Psi$, the map $[\cdot\,,\cdot]^{\mc{G}\,ev,od}$ is $Spin(n,0)\times Spin(0,n)$--equivariant.

\bigskip

Finally, we shall make the following observations for later use. Let $\widehat{m}=(-1)^{m(m+1)/2}$ and $m^{\vee}=(-1)^{m(m-1)/2}$. The Riemannian volume form $\varpi_g$ acts on spinor fields 
\begin{equation}\nonumber
\begin{array}{lcl}
\varpi_g\cdot\Psi_{\pm} & = & \pm\widehat{m}i^m\Psi_{\pm}\\[5pt] \varpi_g\cdot\Psi & = & \phantom{\pm}\widehat{m}i^{m+1}\Psi
\end{array}
,\quad m=[\frac{n}{2}],
\end{equation}
where $\Psi_{\pm}\in\Gamma(\Delta_{\pm})$ are chiral. Since $\star_g\alpha\cdot\Psi=\widehat{\alpha}\cdot\varpi_g\cdot\Psi$ for any form $\alpha$, it follows that
\begin{equation}\label{selfdualgen}
\begin{array}{lcl}
\gtilde[\Psi\otimes\Phi_{\pm}]^{\mc{G}} & = & \pm m^{\vee}i^m[\Psi\otimes\Phi_{\pm}]^{\mc{G}}\\[5pt]\gtilde[\Psi\otimes\Phi]^{\mc{G}} & = &\phantom{\pm}m^{\vee}i^{m+1}\widetilde{[\Psi\otimes\Phi]^{\mc{G}}}
\end{array}
,\quad m=[\frac{n}{2}].
\end{equation}
%
%
%
%
%
\section{Twisted Calabi--Yau structures}\label{gencom}
%
%
%
\subsection{Classical $SU(m)$--structures}
In terms of structure groups, a usual Calabi--Yau structure on $M^{2m}$ is equivalent to a reduction to $SU(m)$ acting in its standard vector representation. We therefore also speak of an $SU(m)$--{\em structure}. Let us first review this reduction mechanism starting from the frame bundle associated with $GL(2m)$. The reduction is achieved by an invariant (almost) symplectic form $\omega\in\Omega^2(M)$ and a complex $m$--form $\Omega\in\Omega^m(M)\otimes\C$, which are subject to various compatibility conditions~\cite{hi97}. In particular, $\Omega$ gives rise to an almost complex structure $J$ inducing a K\"ahler metric $g$. This metric admits a canonical spin structure, for we can lift the standard embedding $SU(m)\hookrightarrow SO(2m)$ to $Spin(2m)$ as $SU(m)$ is simply--connected. 

\bigskip

Next we describe the reduction from $Spin(2m)$ to this lift (cf. for instance Chapter IV.9 in~\cite{lami89}). Recall that a spinor $\Psi\in\Delta$ is said to be {\em pure} if the annihilator $W_{\Psi}=\{z\in\C^{2m}=\R^{2m}\otimes\C\,|\,z\cdot\Psi=0\}$ is maximally isotropic \wrt the complex--linear extension of the inner product on $\R^{2m}$. Pure spinors are necessarily chiral, and we obtain a complex structure corresponding to the decomposition $\R^{2m}\otimes\C=W_{\psi}\otimes\overline{W}_{\psi}$. In fact, there is an $Spin(2m)$--equivariant correspondence between lines of pure spinors and complex structures on $\R^{2m}$. Furthermore, the orbit of pure spinors is isomorphic with $\R_{\not=0}\times Spin(2m)/SU(m)$, where $\R_{\not=0}$ acts on $\Delta_{\pm}$ by rescaling. Hence an $SU(m)$--structure on $M^{2m}$ gives rise to a pure spinor field which we take to be of unit norm. Conversely, any global pure spinor field induces an $SU(m)$--structure on $M$. We recover the symplectic form from the metric and the induced (almost) complex structure, while the complex volume form is obtained by squaring the normalised spinor field, i.e.\ $[\Psi\otimes\Psi]/\eunorm{\Psi}^2=\Omega$.
%
%
%
\subsection{Generalised Calabi--Yau structure}\label{gcy}
In general, squaring a pure spinor field trivialises the canonical line bundle of the induced complex structure, that is, we have an $SL(m,\C)$--structure. Following~\cite{hi03}, these structures generalise as follows.

\begin{definition}
A(n) {\em (almost) generalised complex structure} on $(M^{2m},H)$ is given by an isometry $\mc{J}\in\End\big(\mb{E}\big)$ such that $\mc{J}^2=-Id$. 
\end{definition}

Put differently, a generalised complex structure is an orthogonal splitting $\mb{E}\otimes\C=\mb{W}\oplus\overline{\mb{W}}$ into {\em isotropic} complex vector subbundles $\mb{W}$ and $\overline{\mb{W}}$. In terms of structure groups, this is equivalent to a reduction from $SO(2m,2m)$ to $U(m,m)$. Furthermore, $\mb{W}$ determines the associated line bundle $\mb{K}_{\mc{J}}$ of pure spinor fields in either $\mb{S}_+\otimes\C$ or $\mb{S}_-\otimes\C$. 

\begin{definition}
A(n) {\em (almost) generalised Calabi--Yau structure} on $(M,H)$ is defined by a pure spinor field $\rho\in\Gamma(\mb{S}_{\pm}\otimes\C)$ with $\langle\rho,\overline{\rho}\rangle\not=0$.
\end{definition}

\begin{rmk}
The condition $\langle\rho,\overline{\rho}\rangle\not=0$ ensures that $\mb{W}_{\rho}$ defines a generalised complex structure. Since $\mb{K}_{\mc{J}}\otimes\mb{K}_{\mc{J}}\cong\Lambda^{2m}\mb{W}^*_{\rho}$ for orientable $M$, $\rho$ induces a further reduction from $U(m,m)$ to $SU(m,m)$ (see~\cite{hi03} for details).  
\end{rmk}

\begin{ex}
Consider a classical Calabi--Yau structure $(M,\omega,\Omega)$ whose complex structure is given by $J$. We define the spinor fields $\rho_0=\widehat{m}\exp(-i\omega)$ and $\rho_1=\Omega$, where the sign of $\rho_0$ is put for convenience and has no intrinsic meaning. Then $\rho_0$ and $\rho_1$ trivialise $\mb{K}_{\mc{J}_{\omega}}$ and $\mb{K}_{\mc{J}_J}$ corresponding to 
\begin{equation}\nonumber
\mc{J}_{\omega}=\left(\begin{array}{cc}0&\omega^{-1}\\\omega&0\end{array}\right)\mbox{ and }\,\mc{J}_{\Omega}=\left(\begin{array}{cc} -J & 0\\0 & J^*\end{array}\right).
\end{equation}
Further, the compatibility condition $\Omega\wedge\overline{\Omega}=m^{\vee}i^m\omega^m$ is reflected in the relation 
\begin{equation}\label{length}
\widehat{m}\frac{m!}{2^m}\langle\rho_0,\overline{\rho_0}\rangle=\langle\rho_1,\overline{\rho_1}\rangle\not=0.
\end{equation} 
\end{ex}

In order to generalise $SU(m)$--structures we need an appropriate notion of a generalised K\"ahler metric. The following definitions are due to Gualtieri~\cite{gu04}.

\begin{definition}
{\rm(i)} A(n) {\em (almost) generalised K\"ahler structure} on $(M^{2m},H)$ is defined by a pair of generalised complex structures $(\mc{J}_0,\mc{J}_1)$ such that
\begin{itemize}
	\item $\mc{J}_0\mc{J}_1=\mc{J}_1\mc{J}_0$
	\item $\mc{G}=-\mc{J}_0\mc{J}_1$ defines a generalised metric.
\end{itemize}
\vspace{-5pt}
In terms of structure groups, this corresponds to a reduction from $SO(2m,2m)$ to $U(m)\times U(m)$.

\medskip

{\rm(ii)} A {\em generalised Calabi--Yau metric} on $(M,H)$ consists of $\rho_0\in\Gamma(\mb{S}_+)$, $\rho_1\in\Gamma(\mb{S}_{(-1)^m})$ such that
\begin{itemize}
	\item both $\rho_0$ and $\rho_1$ induce a generalised Calabi--Yau structure with
	\begin{equation}\label{lengthconst}
	c\langle\rho_0,\overline{\rho_0}\rangle=\langle\rho_1,\overline{\rho_1}\rangle
	\end{equation}
	for some constant $c$.
	\item the induced generalised complex structures $(M,\mc{J}_{\rho_0},\mc{J}_{\rho_1})$ define a generalised K\"ahler structure.
\end{itemize}
\end{definition}

\begin{ex}
Taking up our previous example, we see that 
\begin{equation}\nonumber
\mc{G}=-\mc{J}_{\Omega}\mc{J}_{\omega}=-\mc{J}_{\omega}\mc{J}_{\Omega}=\left(\begin{array}{ll}0&g^{-1}\\g&0\end{array}\right).
\end{equation}
Hence, any $SU(m)$--structure defines a generalised Calabi--Yau metric on $M^{2m}$.
\end{ex}
%
%
%
\subsection{Generalised $SU(m)$--structures}\label{gensum}
Generalised Calabi--Yau metrics are defined by a scale--invariant condition. In terms of structure groups of $\mb{E}$, they turn out to be equivalent to the following structure:

\begin{definition} {\rm (cf.~\cite{gmwi06})}
A {\em generalised $SU(m)$--structure} on $(M^{2m},H)$ is a reduction to $SU(m)\times SU(m)\hookrightarrow\R_{\not=0}\times Spin(2m,2m)$, where the inclusion of $SU(m)\times SU(m)$ into $Spin(2m,2m)$ lifts the inclusion $SU(m)\times SU(m)\hookrightarrow SO(2m)\times SO(2m)$ given by the standard vector representation.
\end{definition} 

\begin{lem}
A generalised $SU(m)$--structure is characterised by either set of data:
\begin{itemize}
	\item a metric $g$, a $2$--form $B$, a function $\phi$ and two pure spinor fields $\Psi_{L,R}\in\Gamma(\Delta_+)$ of unit norm. 
	\vspace{-5pt}
	\item an $n$--vector field $\nu$ and a pair of $SU(m)\times SU(m)$--invariant spinor fields $\rho_0\in\Gamma(\mb{S}_+\otimes\C)$, $\rho_1\in\Gamma(\mb{S}_{(-1)^m}\otimes\C)$.
\end{itemize}
If we let $\nu=e^{2\phi}\nu_g$, then
\begin{equation}\nonumber
\rho_0^{\nu}=e^{-\phi}e^B\wedge[\mc{A}(\Psi_L)\otimes\Psi_R],\quad\rho_1^{\nu}=e^{-\phi}e^B\wedge[\Psi_L\otimes\Psi_R]
\end{equation}
for suitably chosen $\Psi_L$ and $\Psi_R$.
\end{lem}

\begin{rmk}
For a given generalised $SU(m)$--structure we usually incorporate the induced $B$--field into the $H$--flux as we already did for generalised metrics (cf. Section~\ref{genriemmet}). For reasons becoming clearer in Section~\ref{rrfields}, we refer to $\phi$ as the {\em dilaton}.
\end{rmk}

Let us briefly outline the argument from~\cite{gmwi06} which is similar to the case of generalised $G_2$--structures~\cite{wi06}: Since the projection of $SU(m)\times SU(m)$ down to $SO(m,m)$ sits inside some $SO(2m)\times SO(2m)\cong SO(\mb{V}^+)\times SO(\mb{V}^-)$, a generalised $SU(m)$--structure singles out a generalised metric $(M,H,g)$. Moreover, the structure groups $SO(\mb{V}^{\pm})\cong SO(2m)$ of the vector bundles $\mb{V}^{\pm}$ reduce to $SU(m)$. Since these are isomorphic with the tangent bundle, the latter is associated with two global $SU(m)$--structures inside the underlying $SO(2m)$--structure, giving rise to the corresponding pure spinor fields $\Psi_{L,R}$. By equivariance of the mapping $[\cdot\,,\cdot]^{\gtilde}$, the $\mb{E}$--spinor fields defined by $\rho_0=[\mc{A}(\Psi_L)\otimes\Psi_R]^{\gtilde}$ and $\rho_1=[\Psi_L\otimes\Psi_R]^{\gtilde}$ are $SU(m)\times SU(m)$--invariant. Furthermore, any such pair $(\rho_0,\rho_1)$ is of this form.

\begin{ex}
Consider an $SU(m)$--structure $(M,g,\Psi)$ given by a pure spinor field $\Psi$ of unit norm. By the lemma, it defines a generalised $SU(m)$--structure with $B=0$, $\phi=0$ and $\Psi_L=\Psi=\Psi_R$, or equivalently
\begin{equation}\nonumber
\rho_0=[\mc{A}(\Psi)\otimes\Psi]=\widehat{m}e^{-i\omega},\quad\rho_1=[\Psi\otimes\Psi]=\Omega.
\end{equation}
\end{ex}

\begin{rmk}
In general, the two $SU(m)$--structures do not reduce to a common subgroup as in the previous example (where they are equal). If they do, the induced generalised $SU(m)$--structure will be referred to as {\em straight}. For instance, consider a Riemannian spin manifold $(M^6,g)$ of dimension $6$ with two pure orthogonal spinor fields $\Psi_L$ and $\Psi_R$ (in fact, in dimension $6$ any non--zero spinor is pure). Then the associated $SU(m)$--structures reduce to $SU(2)$ (the stabiliser of a pair of orthogonal spinors). 
\end{rmk}

\begin{prp}
A generalised Calabi--Yau metric is equivalent to a generalised $SU(m)$--structure.
\end{prp}
\begin{prf}
Let us start with a Calabi--Yau metric $(\rho_0,\rho_1)$. The dilaton accounts for rescaling so we assume the spinor fields to be normalised. Denote by $\mc{J}_0$ and $\mc{J}_1$ the induced endomorphisms of $\mb{E}$. Since these commute, the decomposition $\mb{E}\otimes{\C}=\mb{W}_0\oplus\overline{\mb{W}}_0$ is stable under $\mc{J}_1$. Therefore, restricting $\mc{J}_1$ to $\mb{W}_0$ decomposes this space into $\pm i$--eigenspaces $\mb{W}^{\pm}_0$. As a result, $\mb{E}\otimes\C=\mb{W}^+_0\oplus\mb{W}^-_0\oplus\overline{\mb{W}}^+_0\oplus\overline{\mb{W}}^-_0$. Moreover, since the generalised metric satisfies $\mc{G}=-\mc{J}_0\mc{J}_1$, we obtain an orthogonal decomposition of the definite eigenbundles $\mb{V}^{\pm}\otimes\C=\mb{W}^{\pm}_0\oplus\overline{\mb{W}}^{\pm}_0$. Hence, $\mb{V}^{\pm}$ carry a $U(m)$--structure whose canonical line bundles are $\kappa_{\pm}=\Lambda^m\mb{W}^{\pm *}_0$.
On the other hand, 
\begin{equation}\nonumber
\mb{K}^2_{\mc{J}_0}\cong\Lambda^{2m}\mb{W}_0^*=\Lambda^{2m}(\mb{W}^+_0\oplus\mb{W}^-_0)^*\cong\Lambda^m\mb{W}^{+*}_0\otimes\Lambda^m\mb{W}^{-*}_0=\kappa_+\otimes\kappa_-.
\end{equation}
Since $\mc{J}_1=\mc{G}\mc{J}_0$ acts on $\mb{W}_0^+\oplus\overline{\mb{W}}^-_0$ as $i\cdot Id$, we deduce that
\begin{equation}\nonumber
\mb{K}^2_{\mc{J}_1}\cong\Lambda^{2m}\mb{W}_1^*=\Lambda^{2m}(\mb{W}^+_0\oplus \overline{\mb{W}}_0^-)^*\cong\kappa_+\otimes\overline{\kappa}_-.
\end{equation}
Since $\rho_0$ and $\rho_1$ trivialise $\mb{K}_{\mc{J}_0}$ and $\mb{K}_{\mc{J}_1}$ respectively, $\kappa_+\cong\overline{\kappa}_-\cong\kappa_-$ so that $\kappa_+$ and $\kappa_-$ are trivial. Furthermore,~(\ref{lengthconst}) implies that the normalised trivialisations $\rho_k$ induce constant length trivialisations of $\kappa_{\pm}$, that is, the structure group of $\mb{V}^{\pm}$ reduces to $SU(m)$. Conversely, take an $SU(m)\times SU(m)$--invariant pair $(\rho_0,\rho_1)$. The inclusion $SU(m)\times SU(m)\subset U(m)\times U(m)$ induces a generalised K\"ahler structure $(\mc{J}_0,\mc{J}_1)$ such that the corresponding spaces $\mb{W}_0$ and $\mb{W}_1$ of $(1,0)$--vectors are the annihilators of $\rho_0$ and $\rho_1$. Finally,~(\ref{lengthconst}) follows from~(\ref{length}),~(\ref{commut}) and the $Spin(2m,2m)$--invariance of $\langle\cdot\,,\cdot\rangle$. 
\end{prf}
%
%
%
%
%
\section{Ramond--Ramond fields}\label{rrfields}
In this section we derive the geometric setup which governs supergravity compactifications to a $6$--dimensional space, starting from the $10$--dimensional supersymmetry variations. This involves two kinds of fields, namely the NS-NS-- (Neveu--Schwarz) and R-R-- (Ramond--Ramond) fields. The NS-NS--fields define a generalised $SU(3)$--structure, while the R-R--fields give rise to certain spinor field equations. Before we deal with the mathematical aspects we briefly sketch the content of type II theories, following ~\cite{gsw87}, \cite{luth89}, \cite{po98}, \cite{zw04}.
%
%
%
\subsection{Type II theory}\label{typeiitheo}
Strings are 1--dimensional objects and therefore come in two flavours -- they are {\em open} or {\em closed} if the parametrising intervall is of the form $[a,b]\subset\R$ or $S^1$. As strings evolve in an $n$--dimensional space-time $N^{1,n-1}$, they sweep out a 2--dimensional surface, the {\em worldsheet} $\Sigma$. Instead of thinking of $\Sigma$ as an embedded surface, one can interpret the embedding as a {\em bosonic} field, i.e.\ a map $X_{\Sigma}:\Sigma\to N^{1,n-1}$. Moreover, we also have two {\em fermionic} fields, the so--called {\em left} and {\em right mover} $\psi_L,\,\psi_R\in\Gamma(\Delta^{\Sigma}_{L,R}\otimes TN_{|\Sigma})$, where $\Delta^{\Sigma}_{L,R}$ are spin bundles associated with some spin structure on $\Sigma$ (which do not necessarily coincide). We can use these fields to define an action functional which determines the physical theory. After an appropriate quantisation process, it turns out that for a consistent theory the space--time needs to be of dimension $10$. A key requirement for this action is its invariance under supersymmetry transformations. These map the bosonic field content, which in mathematical terms we can think of as elements in a $Spin(1,9)$--vector representation, in a 1-1--fashion onto the fermionic field content, that is, elements in a $Spin(1,9)$--spinor representation. The amount of supersymmetry one requires to be preserved depends on the theory. For the moment, we know five consistent superstring theories, namely type I, type IIA, type IIB, and two heterotic theories. They are non-trivially related to each other through duality maps.  

\bigskip

In this article we shall consider type II theories. Both types are defined by a {\em closed} string. In this case, $\Sigma$ is diffeomorphic to the cylinder $S^1\times \R$ on which we introduce the complex coordinate $w=\sigma_1+i\sigma_2$. Apart from the trivial spin bundle on $\Sigma$, whose sections are characterised by the {\em Ramond} gluing condition
\begin{equation}\nonumber
\mbox{Ramond } (R):\quad \psi_{L,R}(w+2\pi\sigma_1)=\psi_{L,R}(w),
\end{equation}
there exists a non--trivial bundle with the {\em Neveu--Schwarz} condition
\begin{equation}\nonumber
\mbox{Neveu--Schwarz } (NS):\quad\psi_{L,R}(w+2\pi\sigma_1)=-\psi_{L,R}(w).
\end{equation}
After quantisation the left and right movers become operators over Hilbert spaces which we accordingly label by NS and R. For these Hilbert spaces one can construct a discrete spectrum of states and a mass operator whose eigenvalues assign a ``mass'' to each of these states. 

\bigskip

For a meaningful theory we need a {\em vacuum state}, that is a massless ground state. This makes perfect sense from a phenomenological point of view, as the next level of the mass spectrum has Planck mass. However, particles of our ``real'' world are far too light and should be rather considered as perturbations of the vacuum state. For the action of the zero modes, i.e.\ the operators preserving the massless ground states, a careful analysis shows that the NS--ground states are non--degenerate and can be represented by a vector, while the zero modes acting on the (degenerated) R--ground states satisfy the Clifford algebra relations, reflecting the fact that R--ground states can be represented by a spinor. We single out the physically relevant part of the spectrum by a so--called {\em GSO--projection}. For this we need to introduce the {\em fermion operator} which labels every state with $+$ or $-$. The combination NS--, however, is ruled out as the mass squared of these states is negative.

\bigskip

In type II theories the various particles belong to sectors which are defined by the ordered pair of the zero mode ground states. For instance, the $(NS+,R+)$--sector indicates that the left moving ground state is NS+ while the right moving ground state is R+. To formulate the type II theories consistently we need to impose three selection rules, for instance modular invariance at 1-loop. For type IIA and IIB this forces the vacua to lie in the sectors
\begin{equation}\nonumber
\begin{array}{ll}
\text{IIA:} & (NS+,NS+)\oplus(R+,NS+)\oplus(NS+,R-)\oplus(R+,R-)\\  
\text{IIB:} & (NS+,NS+)\oplus(R+,NS+)\oplus(NS+,R+)\oplus(R+,R+).
\end{array}
\end{equation}
The effective dynamical degrees of freedom for the left and right mover are encapsulated in the normal bundle of the world--sheet in $TN$. Consequently, the internal symmetry group for the degrees of freedom of the NS-- and R--states is $SO(8)$ and $Spin(8)$. By the triality principle we can associate the vector representation $\mb 8_{\text V}$ with NS+ and the two spin representations of positive and negative chirality, $\mb{8}$ and $\mb{8}'$, with R+ and R-- respectively (in the sequel, we shall drop the label $\pm$ for sake of simplicity). The physical fields are then accounted for by elements in the irreducible components of the tensor product associated with a given pair. For instance, the NS-NS--sector is represented by $\mb 8_{\text V}\otimes \mb 8_{\text V}$ and is therefore bosonic. It can be decomposed into the trivial representation~$\mb{1}$, the 2--forms $\Lambda^2\mb{8}_{\rm V}$ and the symmetric 2--tensors $\bigodot^2\mb{8}_{\rm V}$. These modules contain the fluctuations (i.e.\ first--order deformations) of a dilaton field $\phi\in C^{\infty}(N)$, a {\em $B$-field} $B\in\Omega^2(N)$ and a metric $g\in\bigodot^2(N)$ defined on the space--time $N$. The second bosonic sector is the R-R--sector. Here, the elements of the irreducible components induce the so--called {\em Ramond--Ramond potentials}. For type IIA one obtains $C^1\in\Omega^1(N)$ and $C^3\in\Omega^3(N)$, while for type IIB we have $C^0\in\Omega^0(N)$, $C^2\in\Omega^2(N)$ and the self--dual 4--form $C^4_+\in\Omega_+^4(N)$. In the fermionic sectors R-NS and R-NS we find the gravitino $\Psi_{X}$ (which is of spin $3/2$, $X$ denoting the vector index) and the dilatino $\lambda$ (of spin $1/2$). The full description is given in Table~\ref{reps}. The subscript of the gravitino distinguishes with which mover the R--sector is associated; ``$L$'' denotes the left and ``$R$'' the right mover.
\begin{table}[ht]
\begin{center}
\begin{equation}\nonumber
\begin{array}{llcc} 
\text{sector} & \text{representation} & \text{bosons/fermions} & \text{dim}\\
(NS+,NS+) & {\bf 8}_\text{V}\otimes {\bf 8}_\text{V}=\mb{1}\oplus\Lambda^2\mb{8}_{\rm V}\oplus\bigodot^2\mb{8}_{\rm V} & \phi \oplus B \oplus g\mbox{ bosonic } & 1+28+35\\
(R+,NS+) & {\bf 8}\otimes {\bf 8}_\text{V}=\mb{8}'\oplus\Lambda^3\mb{8}' & \lambda_+\oplus \Psi_{X,L}\mbox{ fermionic } & 8+56\\
(R-,NS+) & {\bf 8}'\otimes {\bf 8}_\text{V}=\mb{8}\oplus\Lambda^3\mb{8} & \lambda_+\oplus \Psi_{X,L}\mbox{ fermionic } & 8+56\\
(NS+,R+) & {\bf 8}_\text{V}\otimes {\bf 8}=\mb{8}'\oplus\Lambda^3\mb{8}' & \lambda_-\oplus \Psi_{X,R}\mbox{ fermionic } & 8+56\\
(NS+,R-) & {\bf 8}_\text{V}\otimes {\bf 8}'=\mb{8}\oplus\Lambda^3\mb{8} & \lambda_-\oplus \Psi_{X,R}\mbox{ fermionic } & 8+56\\
(R+,R+) & {\bf 8}\otimes {\bf 8}=\mb{1}\oplus\Lambda^2\mb{8}_{\rm V}\oplus\Lambda^4_+\mb{8}_{\rm V} & C^0 \oplus C^2 \oplus C^4_+\mbox{ bosonic } & 1+28+35\\ 
(R+,R-) & {\bf 8}\otimes {\bf 8}'=\mb{8}_{\rm V}\oplus\Lambda^3\mb{8}_{\rm V} & C^1 \oplus C^3\mbox{ bosonic }  & 8+56\\ 
\end{array}
\end{equation}
\caption{The various sectors appearing in type II theories}
\label{reps}
\end{center}
\end{table}

\bigskip

By considering only massless states and passing to the low energy limit of superstring theory, we effectively enter the realm of supergravity defined on a ten--dimensional Lorentzian space--time $(N,g^{1,9})$. Here, space--time supersymmetry can be achieved by introducing a globally well--defined supersymmetry parameter. This parameter is non--physical in as far it is introduced for mathematical convenience only. For type II theories, $\varepsilon$ is built out of two space--time spinors $\varepsilon_{1,2}$, whence the name type II. These are of opposite chirality for type IIA and of equal chirality for type IIB. For the supergravity action to be invariant under the supersymmetry transformation induced by $(\varepsilon_1,\varepsilon_2)$, the supersymmetry variations of the bosonic and fermionic fields have to vanish. It turns out that the supersymmetry variations of the bosonic fields depend only on $(\varepsilon_1,\varepsilon_2)$ and the physical fermions (the gravitino and the dilatino). Conversely, the variations of the fermionic fields depend only on $(\varepsilon_1,\varepsilon_2)$ and the bosonic field content. We make the usual assumption of a {\em vacuum background}, that is, there are no fermionic fields. This implies the vanishing of the bosonic variations, but of course we still need the vanishing of the fermionic variations,
\begin{equation}\label{varferm}
\delta_{(\varepsilon_1,\varepsilon_2)}\Psi_X=0,\quad \delta_{(\varepsilon_1,\varepsilon_2)}\lambda=0.
\end{equation}
Both equations involve the differentials of the $B$--field and the R-R--potentials $C^p$. More precisely, we have the fields
\begin{equation}\nonumber
H=dB,\quad F^{ev,od}=e^{-B}\wedge d(e^B\wedge C^{od,ev})=d_HC^{od,ev},
\end{equation} 
where we consider the R-R--potentials as an even or odd form in $\Omega^{ev,od}(N)$, following the spirit of generalised geometry. This, of course, must not increase the number of degrees of freedom, as reflected by the anti--self--duality condition on the R-R--fields $F^{ev,od}$
\begin{equation}\label{10hodgedual}
F^{ev,od} = -\star_{1,9}\widehat{F^{ev,od}}.
\end{equation}
Here, we take the Hodge star $\star_{1,9}$ with respect to the Lorentzian space--time metric $g^{1,9}$. 
%
%
%
\subsection{Supersymmetry variations}
Next we discuss the precise shape of the fermionic variations~(\ref{varferm}) based upon the so--called {\em democratic formulation} of Bergshoeff et al.~\cite{berg01}. With our conventions (in particular, $X\cdot X=-g(X,X)\mb{1}$ for Clifford multiplication), these read as
\begin{eqnarray}\label{gravitinovar}
\delta_{\varepsilon}\Psi_X 
& = & 
\nabla_X\varepsilon+\frac{1}{4}(X\llcorner H)\cdot\mc{P}\varepsilon-(-1)^{ev,od}\frac{1}{16}e^{\phi}\sum_{p=ev,od}\widehat{F}^p\cdot X\cdot\mathbb{P}^p\varepsilon,\\
\delta_{\varepsilon}\lambda 
& = & 
d\phi\cdot\varepsilon+\frac{1}{2}H\cdot\mc{P}\varepsilon+\frac{1}{8}e^{\phi}\sum_{p=ev,od}(5-p)\widehat{F}^p\cdot\mathbb{P}^p\varepsilon,\nonumber
\end{eqnarray}
where we rearranged the supersymmetry parameters $\epsilon_{1,2}$ in a $2$--component vector $\varepsilon=(\varepsilon_1,\varepsilon_2)^{\top}$. Further, we use the linear operators $\mc{P},\mathbb{P}^p$ given by
\begin{equation}\nonumber
\mc{P}=\left(\begin{array}{rr} 1 & 0 \\ 0 & -1
\end{array}\right),\,\,\mathbb{P}^p=\left(\begin{array}{rr} 0 & 1
\\ 1 & 0\end{array}\right)\,p\equiv0,\,1\mod4,\,\,\mathbb{P}^p=(-1)^p\left(\begin{array}{rr} 0 & -1 \\ 1 & 0\end{array}\right)\,p\equiv2,\,3\mod4,
\end{equation}
which we consider as elements of $\End(\Delta_+\oplus\Delta_-)$ for type IIA and of
$\End(\Delta_+\oplus\Delta_+)$ for type IIB.

\bigskip

Instead of considering the variations $\delta_{\varepsilon}\Psi_X=0,\,\delta_{\varepsilon}\lambda=0$, we will use the equivalent set of conditions
\begin{equation}\label{variations}
\delta_{\varepsilon}\Psi_X=0,\quad\mu(\delta_{\varepsilon}\Psi)-\delta_{\varepsilon}\lambda=0
\end{equation}
with Clifford multiplication $\mu:\Gamma(TN\otimes\Delta)\to\Gamma(\Delta)$. Thus, the second equation becomes $\sum s_ke_k\cdot \delta_{\varepsilon}\Psi_{e_k}+\delta_{\varepsilon}\lambda=0$ in terms of some fixed local orthonormal basis $e_0,\ldots,\,e_9$ with $s_k=|\!|e_k|\!|$, known as the {\em modified dilatino equation}~\cite{gmpt05}. Using the standard identities $X\cdot F^p=X\wedge F^p-X\llcorner F^p$ and $F^p\cdot X=(-1)^p(X\wedge F^p+X\llcorner F^p)$, we obtain 
\begin{equation}\label{dilatinovarmod}
\mu(\delta_{\varepsilon}\Psi)-\delta_{\varepsilon}\lambda={\rm D}\varepsilon-d\phi\cdot\varepsilon+\frac{1}{4}H\cdot\mc{P}\varepsilon.
\end{equation}
%
%
%
\subsection{Compactification}\label{compact}
{\em Compactification} is a method for constructing space--times $(N,g^{1,9})$ such that equations~(\ref{variations}) hold. For this we take $N$ to be a direct product of the form $(\R^{1,3}\times M^n,g_0^{1,9-n}\times g^n)$, where $g_0^{1,9-n}$ is the flat Minkowski metric on the {\em external} space $\R^{1,9-n}$ and $g=g^n$ some Riemannian metric (to be determined later) on the $n$--dimensional (compact) orientable spin manifold $M$, the {\em internal} space.  With the appropriate ansatz, the conditions for the space--time to be supersymmetric can be entirely phrased in terms of the internal geometry of $M$. We shall discuss the case of type IIB for $n=6$ explicitly (cf. also~\cite{gmpt05},~\cite{masm05}).

\bigskip

The orthogonal decomposition of the tangent bundle gives rise to a factorisation of the corresponding Clifford algebra. More precisely, we have $\cliff(\R^{1,9})\cong\cliff(\R^{1,3})\htimes\cliff(\R^6)$, where $\htimes$ denotes the induced $\Z_2$--graded algebra structure defined on elements of pure degree by $a\htimes b\cdot a'\htimes b'=(-1)^{\deg(b)\deg(a')}a\cdot a'\htimes b\cdot b'$. The spin representation of $\cliff(\R^{1,9})$ is thus isomorphic with the tensor product $\Delta^{1,3}\otimes\Delta^6$, where $\Delta^{1,3}$ and $\Delta^6$ are the spin representations of $\cliff(\R^{1,3})$ and $\cliff(\R^6)$ respectively. Note that the volume forms $\vol^{1,3}$ and $\vol^6$ square to $-\Id$, while $\vol^{1,9}$ squares to $\Id$. By our conventions $\Delta^{1,3}_{\pm}$ are the $\pm i$--eigenspaces of $\vol^{1,3}$ while $\Delta^6_{\pm}$ are the $\mp i$--eigenspaces of $\vol^6$. Hence, we get the decomposition
\begin{equation}\nonumber
\Delta^{1,9}_+=\Delta^{1,3}_+\otimes\Delta^6_+\oplus\Delta^{1,3}_-\otimes\Delta^6_-.
\end{equation}
In passing we remark that that $Spin_0(1,3)\cong SL(2,\C)$, so we have, as for $n=6$, a conjugation map $\mc{A}$ on the space of spinors which interchanges the chirality. For the supersymmetry parameters we choose then the simple splitting
\begin{equation}\label{spinoransatz}
\begin{array}{lcl}
\varepsilon_1 &=& \zeta_+\otimes\Psi_L+\mc{A}(\zeta_+)\otimes\mc{A}(\Psi_L)\\
\varepsilon_2 &=& \zeta_+\otimes\Psi_R + \mc{A}(\zeta_+)\otimes\mc{A}(\Psi_R)\\
\end{array}
\end{equation}
with external and internal unit spinor fields $\zeta_+\in\Gamma(\R^{1,3},\Delta_+^{1,3})$ and $\Psi_{L,R}\in\Gamma(M,\Delta^6_+)$. The former is taken to be parallel with respect to the canonical spin connection induced by the Levi--Civita connection, i.e.\ $\nabla^{g^{1,3}}\zeta_+\equiv0$. This simple ansatz is not only for computational convenience, it also makes sense from a physical point of view as we do not expect to observe supersymmetry directly. Instead, supersymmetry is supposed to be broken. Therefore, one naturally starts from a minimal supersymmetric vacuum background. 

\bigskip

Next we turn to the bosonic field content. We assume the dilaton and the $H$--flux on $N$ to be internal, that is they are induced by $\phi\in C^{\infty}(M)$ and $H\in\Omega^3(M)$. Together with the spinor fields $\Psi_{L,R}$, we therefore get a generalised $SU(3)$--structure for $(M,H)$. For the Ramond--Ramond field $F\in\Omega^{od}(N)$ we introduce two internal fields $F_{a,b}\in\Omega^{od}(M)$. We combine these to preserve $4$--dimensional Poincar\'e invariance, i.e. 
\begin{equation}\label{RRansatz}
F=\vol^{1,3}\wedge F_a+F_b.
\end{equation}
In this setting, the $10$--dimensional Hodge duality constraint (\ref{10hodgedual}) yields on $M$
\begin{equation}\label{6dhodgedual}
F_a=-\gtilde F_b
\end{equation}
whence $\widehat{F}_a=-\gtilde\widehat{F}_b$ as follows from a direct computation or by using Lemma~\ref{formid} in the next section.

\bigskip

Now we can describe the conditions for a supersymmetric type II background in terms of the internal geometry of $M$. Let $\nabla^g$ denote the Levi--Civita connection of $g$ as well as the canonically induced spin connection. We define the metric connections 
\begin{equation}\nonumber
\nabla^{\pm}_XY=\nabla^g_XY\pm\frac{1}{2}H(X,Y,\cdot),
\end{equation} 
whose lift to $\Gamma(\Delta_{\pm})$ we keep on denoting by $\nabla^{\pm}$. 

\begin{prp}\label{intgeo}
{\rm (i)} Any supersymmetric type IIB background compactified to $(M^6,g)$ via the ansatz~(\ref{spinoransatz}) and~(\ref{RRansatz}) is equivalent to a generalised $SU(3)$--structure $(g,\phi,\Psi_L,\Psi_R)$ on $(M,H)$ and an odd form $F_b\in\Omega^{od}(M)$ satisfying the following conditions:
\begin{itemize}
	\item the {\em algebraic constraint}
		\begin{equation}\label{extconst}
		F_b\cdot\Psi_L= 0,\quad\widehat{F}_b\cdot\Psi_R=0.
		\end{equation}
	\item for all $X\in\mf{X}(M)$, the {\em internal gravitino equations}
		\begin{equation}\label{gravconst}
		\begin{array}{lcr}
		\nabla^+_X\Psi_L & = & 
		-\frac{1}{8}e^{\phi}\widehat{F}_b\cdot X\cdot\Psi_R\phantom{.}\\[5pt]
		\nabla^-_X\Psi_R 
		& = &\frac{1}{8}e^{\phi}F_b\cdot X\cdot\Psi_L.
		\end{array}
		\end{equation}
	\item the {\em internal dilatino equations}
		\begin{equation}\label{moddilconst}
		({\rm D}-d\phi+\frac{1}{4}H)\cdot\Psi_L=0,\quad({\rm D}-d\phi-\frac{1}{4}H)\cdot\Psi_R=0.
		\end{equation}	
\end{itemize}

{\rm (ii)} Any supersymmetric type IIA background compactified to $(M^6,g)$ via the ansatz~(\ref{spinoransatz}) (with $\mc{A}(\Psi_R)$ instead of $\Psi_R$ in the definition of $\epsilon_2$) and~(\ref{RRansatz}) is equivalent to a generalised $SU(3)$--structure $(g,\phi,\Psi_L,\Psi_R)$ on $(M,H)$ and an even form $F_b\in\Omega^{ev}(M)$ satisfying the algebraic constraint~(\ref{extconst}), the internal dilatino equation~(\ref{moddilconst}) and for all $X\in\mf{X}(M)$ the gravitino equation
\begin{equation}\nonumber
\begin{array}{lcr}
\nabla^+_X\Psi_L
& = & \frac{1}{8}e^{\phi}F_b\cdot X\cdot\Psi_R\\[5pt]
\nabla^-_X\Psi_R
& = & \frac{1}{8}e^{\phi}\widehat{F}_b\cdot X\cdot\Psi_L 
\end{array}
\end{equation}
\end{prp}
\begin{prf}
(i) Fix coordinates $x^{\mu}$, $\mu=0,\ldots,3$ of $\R^{1,3}$ and $x^k$, $k=4,\ldots,9$ of $M$. We start with the gravitino variation (\ref{gravitinovar}) and consider the ``external'' variation $\delta_{\varepsilon}\Psi_{\partial_{\mu}}$. Since $\zeta_+$ is parallel and $H\in\Omega^3(M)$, we only need to take the R-R part into account. Then, for instance, we find for the Ramond--Ramond $7$--form term
\begin{eqnarray*}
\widehat{F}^7\cdot\delta_{\mu}\cdot\mathbb{P}^3\varepsilon
& = &\phantom{-}(\vol^{1,3}\htimes F^3_a)\cdot(\partial_{\mu}\htimes\mb{1})\cdot\mathbb{P}
^3\varepsilon\\
& = &-(\vol^{1,3}\cdot\partial_{\mu}\htimes F^3_a)\cdot\mathbb{P}
^3\varepsilon\\
& = &\phantom{-}(\partial_{\mu}\cdot\vol^{1,3}\htimes F^3_a)\left(\begin{array}{r}\varepsilon_2\\-\varepsilon_1\end{array}\right).
\end{eqnarray*}
The remaining terms follow analogously. For $\varepsilon_1$ we thus obtain
\begin{eqnarray*}
& &
\phantom{-}\partial_{\mu}\cdot\zeta_+\otimes(F_b^1+F_b^3+F_b^5-iF_a^1-iF_a^3-iF_a^5)\cdot\Psi_L\\
& & 
+\partial_{\mu}\cdot\mc{A}(\zeta_+)\otimes(F_b^1+F_b^3+F_b^5+iF_a^1+iF_a^3+iF_a^5)\cdot\mc{A}(\Psi_L)\\  
& \stackrel{(\ref{6dhodgedual})}{=} &
\partial_{\mu}\cdot\zeta_+\otimes(F_b+i\gtilde F_b)\cdot\Psi_L+\partial_{\mu}\cdot\mc{A}(\zeta_+)\otimes(F_b+i\gtilde F_b)\cdot\mc{A}(\Psi_L).
\end{eqnarray*}
Proceeding in the same way with $\varepsilon_2$, the vanishing of the supersymmetry variation $\delta_{\varepsilon}\Psi_{\mu}=0$ implies the algebraic constraints~(\ref{extconst}), for $\gtilde F_b\cdot\Psi_L=-iF_b\cdot\Psi_L$. The conditions $F_b\cdot\mc{A}(\Psi_L)=0$, $\widehat{F}_b\cdot\mc{A}(\Psi_R)=0$ hold automatically by~(\ref{Aequiv}). 

\bigskip

Next we investigate the gravitino variation $\delta_{\epsilon}\Psi_{\partial_k}$ for internal coordinates. For the NS-NS--part $(\mb{1}\otimes\nabla^g_{\partial_k})\varepsilon+\frac{1}{4}(\partial_k\llcorner H)\cdot\mc{P}\varepsilon$ we find
\begin{equation}\nonumber
\left(\begin{array}{c}\zeta_+\otimes\big(\nabla^g_{\partial_k}+\frac{1}{4}(\partial_k\llcorner H)\big)\cdot\Psi_L+\mc{A}(\zeta_+)\otimes\big(\nabla^g_{\partial_k}+\frac{1}{4}(\partial_k\llcorner H)\big)\cdot\mc{A}(\Psi_L)\\[5pt]
\zeta_+\otimes\big(\nabla^g_{\partial_k}-\frac{1}{4}(\partial_k\llcorner H)\big)\cdot\Psi_R
+\mc{A}(\zeta_+)\otimes\big(\nabla^g_{\partial_k}-\frac{1}{4}(\partial_k\llcorner H)\big)\cdot\mc{A}(\Psi_R)
\end{array}\right).
\end{equation}
The R-R--part can be dealt with in a similar fashion as before. We get
\begin{equation}\nonumber
2\left(\begin{array}{l}\phantom{-}\zeta_+\otimes\widehat{F}_b\cdot\partial_k\cdot\Psi_R
+\mc{A}(\zeta_+)\otimes\widehat{F}_b\cdot\partial_k\cdot\mc{A}(\Psi_R)\\
-\zeta_+\otimes F_b\cdot\partial_k\cdot\Psi_L
-\mc{A}(\zeta_+)\otimes F_b\cdot\partial_k\cdot\mc{A}(\Psi_L)
\end{array}\right).
\end{equation}
Combining the NS-NS-- and R-R--part, we finally derive the internal gravitino equation~(\ref{gravconst}).

\bigskip

For the modified dilatino equation (\ref{dilatinovarmod}) we note that there is no external contribution, for the spinors fields $\zeta_{\pm}$ are parallel and the $H$--field and the dilaton $\phi$ are only defined on the internal space. Thus the modified dilatino variation boils down to~(\ref{moddilconst}).

(ii) Mutatis mutandis, the same procedure can be carried out for type IIA. 
\end{prf}
%
%
%
\subsection{A no--go theorem}
To illustrate the relevance of Ramond--Ramond fields from a mathematical perspective, we investigate the internal gravitino and dilatino equation for the case $F_b=0$.

\begin{ex}
In absence of any bosonic fields, the conditions~(\ref{extconst}),~(\ref{gravconst}) and~(\ref{moddilconst}) reduce to $\nabla^g\Psi_{L,R}\equiv0$. In particular, the angle $q(\Psi_L,\Psi_R)$ is constant. As a result, we obtain a straight generalised $SU(m)$--structure whose holonomy is contained in $SU(3)$ and reduces even further to $SU(2)$ unless $\Psi_L=\pm\Psi_R$.
\end{ex}

\bigskip

We are going to prove that $F_b=0$ prevents the existence of a solution other than a straight one provided $M$ is compact. In fact, the subsequent arguments apply to any $n$--dimensional manifold $M$ with $(H,g,\phi,\Psi)$ such that the spinor field equations
\begin{equation}\label{gravdil}
\nabla^+\Psi=\nabla_X^g\Psi+\frac{1}{4}(X\llcorner H)\cdot\Psi=0,\quad(d\phi+\frac{1}{2}H)\cdot\Psi=0,
\end{equation}
hold. The key assumption here is that $H$ is closed. Instances other than the internal spaces of the previous section are strongly integrable generalised $G_2$-- or $Spin(7)$--manifolds~\cite{wi06}. 

\bigskip

First we need a curvature result. Applying Corollary 5.2 of~\cite{friv02}, the Ricci endomorphism $\ric^+$ of $\nabla^+$ is (using $\nabla^+_X(H\cdot\Psi)=(\nabla^+_XH)\cdot\Psi$ etc.)
\begin{equation}\nonumber
\ric^+(X)\cdot\Psi=(\nabla^+_XH)\cdot\Psi=-2(\nabla_X^+d\phi)\cdot\Psi,
\end{equation}
whence $\ric^+(X)=-2\nabla_X^+d\phi$. 

\begin{lem}\label{scal}
The scalar curvature of $\nabla^+$ is $S^+=2\Delta^g\phi$, where $\Delta^g$ is the Riemannian Laplacian on functions.
\end{lem}
\begin{prf}
Picking a frame with $\nabla_{e_i}e_j=0$ at a fixed point, we obtain
\begin{eqnarray*}
\ric^+(e_j,e_k) & = &
-2g(\nabla_{e_j}^+d\phi,e_k)\\
& = & -2e_j.g(d\phi,e_k)+g(d\phi,\nabla^+_{e_j}e_k)\\
& = & -2e_j.e_k.\phi+(e_k\llcorner e_j\llcorner H).\phi/2.
\end{eqnarray*}
The first summand is minus twice $\mc{H}^{\phi}$, the Hessian of $\phi$ evaluated in the basis $\{e_k\}$. Hence $\ric^+(X,Y)=-2\mc{H}^{\phi}(X,Y)-X\llcorner Y\llcorner H/2$, so that $S^+=2\Delta^g\phi$.
\end{prf}

\begin{prp}\label{dilH}
For a spinor field $\Psi\in\Gamma(\Delta)$ satisfying~(\ref{gravdil}), it follows $S^+=-3\!\eunorm{H}^2$. In particular, $\phi=0$ \iff $H=0$. Furthermore, if $M$ is compact, then $H=0$.
\end{prp}
\begin{prf}
Let $\D^+$ denote the Dirac operator associated with $\nabla^+$, i.e.\ locally $D^+\Psi=\sum e_k\cdot\nabla^+_{e_k}\Psi$. By Theorem 3.3 in~\cite{friv02},
\begin{equation}\nonumber
\D^+(H\cdot\Psi)=\big(d^*\!H\cdot-2\sigma^H\cdot-2\sum(e_k\llcorner H)\cdot\nabla^+_{e_k}\big)\Psi,
\end{equation}
where  $2\sigma^H=\sum(e_k\llcorner H)\wedge(e_k\llcorner H)$. Now
\begin{eqnarray*}
-2\sum (e_k\llcorner H)\cdot\nabla^+_{e_k}\cdot\Psi & = & \frac{1}{2}\sum(e_k\llcorner H)\cdot(e_k\llcorner H)\cdot\Psi\\
& = &(\frac{3}{2}\eunorm{H}^2+\sigma^H+\alpha)\cdot\Psi
\end{eqnarray*}
for some further $2$--form $\alpha$. On the other hand side, using the proof of the previous lemma,
\begin{eqnarray*}
D^+(H\cdot\Psi) & = & -2\sum e_k\cdot(\nabla^+_{e_k}d\phi)\cdot\Psi\\
& = & -2\sum (e_k\wedge\nabla^+_{e_k}d\phi-e_k\llcorner\nabla^+_{e_k}d\phi)\cdot\Psi\\
& = & (\alpha'-2\Delta^g\phi)\cdot\Psi,
\end{eqnarray*}
for some $2$--form $\alpha'$. Hence contracting with $q(\Psi,\cdot)$ yields
\begin{equation}\nonumber
q\big(\Psi,-2(\Delta^g\phi)\cdot\Psi\big)=q\big(\Psi,(\frac{3}{2}\!\eunorm{H}^2-\sigma^H)\cdot\Psi\big)
\end{equation}
\big($q(\alpha^p\cdot\Psi,\Psi)$ is purely imaginary for $p\equiv2(4)$ and real for $p\equiv0(4)$\big). As $q(\Psi,2\sigma^H\cdot\Psi)=S^+$ by Corollary 3.2 in~\cite{friv02}, the previous lemma implies $3\eunorm{H}^2=-S^+$. Hence, if $M$ is compact, $H=0$, as follows from integration of $\eunorm{H}^2$. 
\end{prf}

\begin{rmk}
There are compact examples of generalised $G_2$-- and $Spin(7)$--structures satisfying~(\ref{gravdil}) for $H,\,dH\not=0$~\cite{wi06}.
\end{rmk}
%
%
%
%
%
\section{Variational principles and integrability}\label{varprinc}
%
%
%
\subsection{Generalised Hodge theory}
Let $(M^n,\nu,H,g)$ be a closed oriented generalised Riemannian manifold. Integration of $\mc{Q}^{\pm}$ (cf. Section~\ref{gengeo}) with respect to the volume form $\varpi_{\nu}$ dual to $\nu$ yields an inner product on $\Gamma(\mb{S}_{\pm})$ given by
\begin{equation}\nonumber
(\sigma,\tau)=\frac{1}{2}\int_M\mc{Q}^{\pm}(\sigma,\tau)\varpi_{\nu}=\frac{1}{2}\int_M[\star\sigma^{\nu}\wedge\tau^{\nu}]^n.
\end{equation}
As we have remarked above, the complex $d_{\nu}:\Gamma(\mb{S}_{\pm}\to\Gamma(\mb{S}_{\mp})$ is elliptic so that Hodge theory applies. Using the straightforward

\begin{lem}\label{formid}
For a differential $p$--form $\alpha$ and a vector field $X$ (respectively its dual), we have the identities
\begin{equation}\nonumber
\star\widehat{\alpha}=(-1)^{[\frac{n}{2}]-np}\widehat{\widetilde{\star\alpha}}=\left\{\begin{array}{ll}(-1)^{m+p}\widehat{\star\alpha},&n=2m\\(-1)^{m+1}\widehat{\star\alpha},&n=2m+1\end{array}\right.\!\!\!\!,\quad\star(X\wedge\alpha)=X\llcorner\star\widetilde{\alpha},\quad\widehat{d\alpha}=-d\widehat{\widetilde{\alpha}}.
\end{equation} 
\end{lem}

one immediately shows that the formal adjoint $d^*_{\nu}$ of $d_{\nu}$ is
\begin{equation}\nonumber
d^*_{\nu}=(-1)^{[\frac{n}{2}]+n}\gtilde d_{\nu}\gtilde.
\end{equation}

\begin{rmk}
In particular, $d_H^*=\gtilde d_H\gtilde$ is the formal adjoint of $d_H$ \wrt the standard $L^2$--inner product on forms. In absence of $H$ this gives the usual formal adjoint $d^*=(-1)^{np+n+1}\!\star d\,\star$ for forms of pure degree $p$.
\end{rmk}

Hence, any cohomology class in $H^{\pm}(\mb{E})$ has a unique harmonic representative $\tau\in\Gamma(\mb{S}_{\pm})$, i.e.\ $d_{\nu}\tau=0$, $d_{\nu}\gtilde\tau=0$. In analogy to classical Hodge theory, these harmonic spinor fields can be thought of as critical points of the functional $\mb{Q}(\tau)=(\tau,\tau)$ restricted to a cohomology class.

\bigskip

In order to incorporate a Ramond--Ramond field, we associate with $\gamma\in\Gamma(\mb{S}_{\pm})$ the functional $\tau\in\Gamma(\mb{S}_{\pm})\mapsto\mb{C}_{\gamma}(\tau)=(\gtilde\gamma,\tau)$. Restricted to the cohomology class of $\tau$, the first variation of $\mb{Q}$ at $\tau$ is $\delta\mb{Q}(\dot{\tau})=(\tau,\dot{\tau})$, where $\dot{\tau}=d_{\nu}\sigma$ is a spinor tangent to $[\tau]$. On the other hand, $\delta\mb{C}(\dot{\tau})=(\gtilde\gamma,\dot{\tau})$ so that introducing a Lagrange multiplier yields

\begin{prp}
A $d_{\nu}$--closed chiral $\mb{E}$--spinor field $\tau$ is a critical point in its cohomology class for $\mb{Q}$ subject to $\mb{C}_{\gamma}(\tau)=const$ if and only if $d_{\nu}\gtilde\tau=\lambda d_{\nu}\gamma$, $\lambda\in\R$. 
\end{prp}

\bigskip

The variational principle we have discussed so far presupposes the choice of a generalised metric. For $n=6$, the case relevant for physics, there exists an intrinsic variational problem for generalised Calabi--Yau structures relying on Hitchin's notion of stability~\cite{hi03}. Recall that such a structure is specified by a pure complex spinor field $\rho=\tau+i\tau^{\dagger}\in\Gamma(\mb{S}_{\pm}\otimes\C)$ (cf. Section~\ref{gcy}) with $\langle\rho,\overline{\rho}\rangle\not=0$. The real part of $\rho$ is a section of the fibre bundle associated with the {\em open} orbit $\R_{>0}\times Spin(6,6)_+/SU(3,3)$. Consequently, the set $\mc{U}$ of real spinor fields $\tau$ defining a generalised Calabi--Yau structure is an open submanifold of $\Gamma(\mb{S}_{\pm}\otimes\C)$. Following Hitchin's language such spinor fields are called {\em stable}. A stable $\tau\in\mc{U}$ determines $\tau^{\dagger}\in\mc{U}$ in a non--linear fashion. Note that $\tau^{\dagger\dagger}=-\tau$. For instance, if $\rho$ happens to be the component of a generalised $SU(3)$--structure, then $\tau^{\dagger}=\gtilde\tau$, where $\gtilde$ is induced by the underlying generalised $SU(3)$--structure). We define a functional $\mb{V}:\mc{U}\to\R$ by
\begin{equation}\nonumber
\mb{V}(\tau)=\int_M\langle\tau^{\dagger},\tau\rangle\varpi_{\nu}=\frac{1}{2i}\int_M\langle\rho,\overline{\rho}\rangle\varpi_{\nu}.
\end{equation}
Since stability is an open condition, we can differentiate this functional and ask for a critical points in a given a cohomology class. In~\cite{hi03}, Hitchin shows that $(\delta V)_{\tau}(\dot{\tau})=\int_M\langle\tau^{\dagger},\dot{\tau}\rangle$,  so that a $d_{\nu}$--closed spinor field $\tau$ is a critical point in its cohomology class \iff $d_{\nu}\tau^{\dagger}=0$. In particular, we obtain $d_{\nu}\rho=0$, which Hitchin adopts as an integrability condition for generalised Calabi--Yau structures. If a generalised $SU(3)$--structure $(\rho_0,\rho_1)$ consists of two such integrable Calabi--Yau structures, then $d\rho_{0,1}=0$, so that $\tau_{0,1}$ is harmonic with respect to the induced generalised metric. 

\bigskip

Again we can constrain the variational principle. For $\gamma\in\Gamma(\mb{S}_{\pm})$ we define the functional $\mb{C}_{\gamma}(\tau)=\langle\gamma,\tau\rangle$ on $\Gamma(\mb{S}_{\pm})$. As above, we deduce the

\begin{prp}
A closed chiral spinor field $\tau\in\Gamma(\mb{S}_{\pm})$ is a critical point in its cohomology class for $\mb{V}$ subject to $\mb{C}_{\gamma}(\tau)=const$ \iff $d_{\nu}\tau^{\dagger}=\lambda d_{\nu}\gamma$, $\lambda\in\R$. 
\end{prp}

\begin{ex}
Consider a straight generalised $SU(3)$--structure given by $\rho_0=e^{-i\omega}=1-i\omega-\omega^2/2+i\omega^3/6$ and $\rho_1=\Omega=\psi_++i\psi_-$ as in the example of Section~\ref{gensum}. Then both $\rho_0$ and $\rho_1$ are induced by an unconstrained critical point of $\mb{V}$ \iff $d\omega=0$, $d\psi_{\pm}=0$. Equivalently, the holonomy of $g$ is contained in $SU(3)$ (cf. for instance~\cite{hi97}). Further, if $d\omega\not=0$, $\tau_0=1-\omega^2/2$ defines a constrained critical point for any constant multiple of $\gamma=\omega+d\alpha^1$ provided $d\omega\wedge\omega=0$. Similarly, if $d\psi_{\mp}\not=0$, then $\tau_1=\psi_{\pm}$ defines a constrained critical point for any constant multiple of $\gamma=\psi_{\mp}+d\alpha^2$ provided $d\psi_{\pm}=0$. We come back to this example in the next section.
\end{ex}
%
%
%
\subsection{Spinor field equations}
In the light of the constrained variational principle, we consider a generalised $SU(m)$--structure $(\rho_0,\rho_1)$ on $(M,H)$ ($M$ not necessarily compact) satisfying $d_{\nu}\rho_{0,1}=\varphi_{0,1}$ for (not necessarily exact) complex--valued spinor fields $\varphi_{0,1}$. We can rephrase this condition in terms of the corresponding data $(g,\phi,\Psi_L,\Psi_R)$ (as before, we modify $H$ by the global $B$--field coming from the generalised $SU(m)$--structure if necessary). Let $F_{0,1}=\varphi_{0,1}^{\nu}$, then
\begin{equation}\label{intcondrr}
d_He^{-\phi}[\mc{A}(\Psi_L)\otimes\Psi_R]=F_0,\quad d_He^{-\phi}[\Psi_L\otimes\Psi_R]=F_1.
\end{equation}

\bigskip

Before we try to reformulate these equations in terms of spinor field equations on $\Psi_{L,R}$, we put an additional condition on the forms $F_{0,1}$ which for $n=6$ will reproduce the constraint~(\ref{extconst}). To this end, assume $m\geq3$ and let $\lambda^k=\Lambda^k\C^m$ be the $k$--th exterior power of the vector representation $\C^m$ of $SU(m)$. Then $\lambda^m\cong\C$ and $\overline{\lambda^{m-k}}\cong(\lambda^{m-k})^*\cong\lambda^k$ . Further, we have an orthogonal decomposition into $SU(m)$--representation spaces  $\Delta_+=\lambda^m\oplus\lambda^{m-2}\oplus\ldots$ and $\Delta_-=\lambda^{m-1}\oplus\lambda^{m-3}\oplus\ldots$ (cf. for instance~\cite{wa89}). Schematically\begin{eqnarray*}
m\mbox{ even: }&\Delta_+=\C\oplus W\oplus\C,&\Delta_-=\overline{\C^m}\oplus \widetilde{W}\oplus\C^m,\\
m\mbox{ odd: }&\Delta_+=\C\oplus V\oplus\C^m,&\Delta_-=\overline{\C^m}\oplus\overline{V}\oplus\C,
\end{eqnarray*}
where $V$, $W$ and $\widetilde{W}$ are (not necessarily irreducible) $SU(m)$--modules. 

\begin{ex}
For $m=3$ we find $V=\{0\}$, while for $m=4$, $W=\lambda^2$ and $\widetilde{W}=\{0\}$. 
\end{ex}

\bigskip

The $SU(m)$--invariant spinors spanning the trivial representations are $\Psi$ and $\mc{A}(\Psi)$. As noted in Section~\ref{gencom}, $\Psi$ \big(and thus $\mc{A}(\Psi)$\big) are pure. The induced complex structure $J$ satisfies
\begin{equation}\label{Jaction}
X\cdot\Psi=iJX\cdot\Psi,
\end{equation} 
$X\in\R^{2m}$. Taking $SU(m)=SU(m)_L$ with invariant spinors $\Psi_L$ and $\mc{A}(\Psi_L)$, an orthonormal basis of $\C^m$ (resp.\ $\overline{\C^m}$) is given by $z_{L,k}\cdot\mc{A}(\Psi_L)$ (resp. $\overline{z}_{L,k}\cdot\Psi_L$) for a unitary basis $z_{L,k}=(e_k-iJ_Le_k)/2$ of $\C^m$. Replacing $\Psi_L$ by $\Psi_R$ and $J_L$ by $J_R$ gives the decomposition as an $SU(m)_R$--module. Further, 
\begin{equation}\nonumber
\begin{array}{l}
\Lambda^{ev} = \Delta_+\otimes\Delta_+\oplus\Delta_-\otimes\Delta_-,\,\Lambda^{od} = \Delta_+\otimes\Delta_-\oplus\Delta_-\otimes\Delta_+\mbox{ for $m$ even}\\[5pt]
\Lambda^{ev}=\Delta_+\otimes\Delta_-\oplus\Delta_-\otimes\Delta_+,\,\Lambda^{od}=\Delta_+\otimes\Delta_+\oplus\Delta_-\otimes\Delta_-\mbox{ for $m$ odd}. 
\end{array}
\end{equation}
In the sequel, given $SU(m)_L\times SU(m)_R\cong SU(m)\times SU(m)$, we regard the left resp. right hand side factors of $\Delta\otimes\Delta$ as an $SU(m)_L$-- resp. $SU(m)_R$--module. We shall write $\C^m_L$, $\C^m_R$ etc.\ accordingly. In particular, all these decompositions acquire global meaning for $(M,H)$ endowed with a generalised $SU(m)$--structure.

\begin{definition}
Consider a generalised $SU(m)$--structure given by $(g,\phi,\Psi_L,\Psi_R)$.

{\rm (i)} Let $m$ be even. An odd form $F\in\Omega^{od}(M)\otimes\C$ will be called a {\em Ramond--Ramond field} for the generalised $SU(m)$--structure if pointwise
\begin{equation}
F\in\widetilde{W_L}\otimes\overline{\C^m_R}\oplus\overline{\C^m_L}\otimes W_R\subset\Delta_+\otimes\Delta_-\oplus\Delta_-\otimes\Delta_+.
\end{equation}

{\rm (ii)} Let $m$ be odd. An odd form $F\in\Omega^{od}(M)\otimes\C$ will be called a {\em Ramond--Ramond field} for the generalised $SU(m)$--structure if pointwise
\begin{equation}
F\in\C^m_L\otimes\big(\C^m_R\oplus V_R\big)\oplus\big(\overline{\C^m_L}\oplus\overline{V_L}\big)\otimes\overline{\C_R^m}\subset\Delta_+\otimes\Delta_+\oplus\Delta_-\otimes\Delta_-.
\end{equation}
Similarly, an even complex form $F\in\Omega^{ev}(M)\otimes\C$ will be a Ramond--Ramond field if pointwise
\begin{equation}\nonumber
F\in\big(\C^{m}_L\oplus V_L\big)\otimes\overline{\C^m_R}\oplus\overline{\C^m_L}\otimes\big(\C^m_R\oplus V_R\big)\subset\Delta_+\otimes\Delta_-\oplus\Delta_-\otimes\Delta_+.
\end{equation}
\end{definition}

\bigskip

The following proposition justifies our jargon.
 
\begin{prp}\label{rrisrr}
For any Ramond--Ramond field associated with a generalised $SU(m)$--structure $(\rho_0,\rho_1)\cong(g,\phi,\Psi_R,\Psi_L)$ on $(M,H)$ we have
\begin{equation}\label{rrchar}
\widehat{F}\cdot\Psi_L=0,\;\widehat{F}\cdot\mc{A}(\Psi_L)=0,\quad F\cdot\Psi_R=0,\;F\cdot\mc{A}(\Psi_R)=0.
\end{equation}
Equivalently,
\begin{equation}\nonumber
\widehat{F}^{\pm}\bullet\rho_{0,1}=0,
\end{equation}
where $F^{\pm}$ denotes the image of $F$ under the extension of $X\in TM\mapsto X^{\pm}\in\mb{V}^{\pm}$ (cf. Section~\ref{genriemmet}). Furthermore, the converse is true for $m=3$. 
\end{prp}
\begin{prf}
First, we need a technical lemma.

\begin{lem}
Let  $d=\dim_{\C}\Delta=2^m$ and $\Psi,\,\Phi,\,\xi\in\Delta$. Then
\begin{equation}\label{faction}
\frac{1}{d}\widehat{[\Psi\otimes\Phi]}\cdot\xi=q(\mc{A}(\Psi),\xi)\Phi.
\end{equation} 
and
\begin{equation}\label{hatbar}
\widehat{[\Phi\otimes\Psi]}=\left\{\begin{array}{ll}m\mbox{ even:}\quad&\widehat{m}\widetilde{[\Psi\otimes\Phi]}\\[5pt]
m\mbox{ odd:}\quad&\widehat{m}[\Psi\otimes\Phi]\end{array}\right..
\end{equation}
\end{lem}
\begin{prf}
The trace operator ${\rm Tr}(A^*B)/d$ defines a positive--definite hermitian inner product on $\End(\Delta)$. If $\{e_K=e_{k_1}\cdot\ldots\cdot e_{k_r}\}$ is an orthonormal basis for $\cliff(W,g)$ and thus for $\cliff(W^{\C},g^{\C})$, then so is $\{\kappa(e_K)\}$ for $(\End(\Delta),{\rm Tr}/d)$. First we note that $\kappa(a)^*=\kappa(\widehat{\overline{a}})$ and $g^{\C}(a,b)=[\widehat{a}\cdot b]^0$, that is, $g^{\C}(a,b)$ equals the zero degree component of $\widehat{a}\cdot b$. Let $c$ such that $g^{\C}(c,a)={\rm Tr}\big(\kappa(a)\big)$. Then $c$ belongs to the centre of $\cliff(W^{\C},g^{\C})$ (for $g^{\C}(a\cdot b,c)=g^{\C}(b,\widehat{a}\cdot c)=g^{\C}(a,c\cdot \widehat{b})$), whence $c=d$. In particular,  ${\rm Tr}\big(\kappa(e_K)^* \kappa(e_L)\big)/d=g^{\C}(1,\widehat{e_K}\cdot e_L)=\delta_{KL}$. Next let $\{\xi_k\}$ be an orthonormal basis of $\Delta$. Define $\Psi\odot\Phi:\Delta\to\Delta$ by $\Psi\odot\Phi(\xi)=q(\mc{A}(\Psi),\xi)\Phi$. Then
\begin{eqnarray*}
\Psi\odot\Phi
&=&\frac{1}{d}\sum_K{\rm Tr}\big(\kappa(e_K)^*\circ\Psi\odot\Phi\big)\kappa(e_K)\\
&=&\frac{1}{d}\sum_{K,k}q\big(\mc{A}(\Psi),\xi_k\big)q\big(e_K\cdot\xi_k,\Phi\big)\kappa(e_K)\\
&=&\frac{1}{d}\sum_Kq\big(\mc{A}(\Psi),e_K\cdot\Phi\big)\kappa(\widehat{e_K})\\
&=&\frac{1}{d}\widehat{[\Psi\otimes\Phi]},
\end{eqnarray*}
which proves~(\ref{faction}). Further,~(\ref{hatbar}) is a straightforward consequence of $q\big(\mc{A}(\Psi),\mc{A}(\Phi)\big)=\overline{q(\Psi,\Phi)}$.
\end{prf}

\begin{cor}\label{contlemma}
{\rm (i)} Contraction of $\Psi\otimes\Phi$ with $q(\xi,\cdot)\otimes\mathrm{Id}$ yields
\begin{equation}\nonumber
q(\xi,\Psi)\Phi=q\big(\mc{A}(\Psi),\mc{A}(\xi)\big)\Phi=\frac{1}{2^m}\widehat{[\Psi\otimes\Phi]}\cdot\mc{A}(\xi).
\end{equation}

{\rm (ii)} Contraction of $\Psi\otimes\Phi$ with $\mathrm{Id}\otimes q(\xi,\cdot)$ yields
\begin{equation}\nonumber
q(\xi,\Phi)\Psi=\frac{1}{2^m}\widehat{[\Phi\otimes\Psi]}\cdot\mc{A}(\xi)=\left\{\begin{array}{ll}m\mbox{ even:}\quad&\frac{1}{2^m}\widehat{m}\widetilde{[\Psi\otimes\Phi]}\cdot\mc{A}(\xi)\\[5pt]
m\mbox{ odd:}\quad&\frac{1}{2^m}\widehat{m}[\Psi\otimes\Phi]\cdot\mc{A}(\xi)\end{array}\right..
\end{equation}
\end{cor}

To be concrete take $m$ to be odd and $F\in\Omega^{od}(M)$. Contraction of $F$ with $\mathrm{Id}\otimes q(\Psi_R,\cdot)$ kills off all components since the right hand side factors are orthogonal to $\Psi_R$. By the previous corollary, the contraction is equal to $F\cdot\mc{A}(\Psi_R)$ up to some constant and therefore vanishes. The remaining cases follow similarly. Further, this is equivalent to $F^{\pm}\bullet\rho_{0,1}=0$ by~(\ref{commut}). To show the converse for $m=3$, note that contraction of $F$ with $\mathrm{Id}\otimes q(\Psi_R,\cdot)$ kills off all components of $F$ except the ones in $\C_L\otimes\C_R$ and $\C^3_L\otimes\C_R$ which get mapped to $\C_L$ and $\C^3_L$. Hence $F\cdot\mc{A}(\Psi_R)=0$ forces these components to vanish. The other cases follow in the same fashion.
\end{prf}

\begin{rmk}
In the light of this proposition, the integrability condition~(\ref{intcondrr}) with Ramond--Ramond fields $F_0$ and $F_1$ is opposite to the one defining so--called {\em weakly integrable} generalised $G_2$--structures~\cite{wi06}, where $F$ is only allowed to have components in the trivial representations spanned by $\Psi_L\otimes\Psi_R$ etc..
\end{rmk}

\bigskip

We now state our central result. Consider the twisted Dirac operators on $\Gamma(\Delta\otimes\Delta)$ given by
\begin{eqnarray*}
\mc{D}(\Psi_1\otimes\Psi_2) & = & \sum
e_k\cdot\nabla^g_{e_k}\Psi_1\otimes\Psi_2+e_k\cdot\Psi_1\otimes\nabla^g_{e_k}\Psi_2\\
& = & \D\Psi_1\otimes\Psi_2+\sum e_k\cdot\Psi_1\otimes\nabla^g_{e_k}\Psi_2,\\
\widetilde{\mc{D}}(\Psi_1\otimes\Psi_2) & = & \sum \nabla^g_{e_k}\Psi_1\otimes
e_k\cdot\Psi_2+\Psi_1\otimes e_k\cdot\nabla^g_{e_k}\Psi_2\\ 
& = & \sum\nabla^g_{e_k}\Psi_1\otimes e_k\cdot\Psi_2+\Psi_1\otimes\D\Psi_2.
\end{eqnarray*}
for any local orthonormal basis $e_1,\ldots,e_{2m}$. As above, $\D:\Gamma(\Delta_{\pm})\to\Gamma(\Delta_{\mp})$ is the natural Dirac operator associated with the spin structure. 

\begin{thm}\label{integrability}
Let $(g,\,\phi,\,\Psi_L,\,\Psi_R)$ be a generalised $SU(m)$--structure with Ramond--Ramond fields  $F_0$ and $F_1$. Then
\begin{equation}\nonumber
d_He^{-\phi}[\mc{A}(\Psi_L)\otimes\Psi_R]=F_0,\quad d_He^{-\phi}[\Psi_L\otimes\Psi_R]=F_1
\end{equation}
holds \iff the {\em gravitino equations}
\begin{eqnarray*}
\nabla_X^+\Psi_L & = & \phantom{-}\frac{1}{2^m}e^{\phi}\big(\widehat{m}\overline{F_0}\cdot X\cdot\Psi_R+m^{\vee}F_1\cdot X\cdot\mc{A}(\Psi_R)\big)\nonumber\\
\nabla_X^-\Psi_R & = & -\frac{1}{2^m}e^{\phi}\big(\widehat{m}\widehat{F_0}\cdot X\cdot\Psi_L+\widehat{F_1}\cdot X\cdot\mc{A}(\Psi_L)\big)
\end{eqnarray*}
and the {\em dilatino equations}
\begin{equation}\nonumber
(\D-d\phi\pm\frac{1}{4}H\big)\cdot\Psi_{L,R}=0
\end{equation}
hold for all $X\in\mf{X}(M)$.
\end{thm}

\begin{cor}
A generalised $SU(m)$--structure $(\rho_0,\rho_1)$ on $(M,\nu,H)$ together with a $d_H$--exact Ramond--Ramond field $F_b$

{\rm (i)} gives rise to a supersymmetric type IIB background compactified to $M$ via the ansatz~(\ref{spinoransatz}) and~(\ref{RRansatz}) \iff $\rho_0/\rho_1$ is a constrained/unconstrained critical point, i.e.\ $d\rho_0=-\widehat{F}_b$ or $d\rho_0=i\widehat{F}_b$ and $d\rho_1=0$.

{\rm (ii)} gives rise to a supersymmetric type IIA background compactified to $M$ via the ansatz~(\ref{spinoransatz}) and~(\ref{RRansatz}) \iff $\rho_0/\rho_1$ is an unconstrained/constrained critical point, i.e.\ $d\rho_0=0$ and $d\rho_1=-F_b$ or $d\rho_1=iF_b$.
\end{cor}

\begin{prf}
Let us assume that the generalised $SU(m)$--structure is integrable. Then
\begin{eqnarray*}
d[\mc{A}(\Psi_L)\otimes\Psi_R] & = & e^{\phi}F_0-(\alpha+H)\wedge [\mc{A}(\Psi_L)\otimes\Psi_R]\\
d[\Psi_L\otimes\Psi_R] & = & e^{\phi}F_1-(\alpha+H)\wedge[\Psi_L\otimes\Psi_R],
\end{eqnarray*}
where we put $\alpha=-d\phi$. We note that $\widetilde{F}_0=-F_0$ and $\widetilde{F}_1=(-1)^{m+1}F_1$. 
By~(\ref{selfdualgen}), $\gtilde[\Psi_L\otimes\Psi_R]=m^{\vee}i^m[\Psi_L\otimes\Psi_R]$, whence
\begin{equation}\nonumber
d\star\widehat{[\Psi_L\otimes\Psi_R]}=m^{\vee}i^m\big(e^{\phi}F_1-(\alpha+H)\wedge[\Psi_L\otimes\Psi_R]\big)
\end{equation}
and similarly for $[\mc{A}(\Psi_L)\otimes\Psi_R]$. By Lemma~\ref{formid}, it follows
\begin{eqnarray*}
d^*[\mc{A}(\Psi_L)\otimes\Psi_R] & = & \widehat{m}i^me^{\phi}\gtilde F_0+(\alpha+H)\llcorner [\mc{A}(\Psi_L)\otimes\Psi_R],\\
d^*[\Psi_L\otimes\Psi_R] & = & \widehat{m}i^me^{\phi}\gtilde F_1+(\alpha+H)\llcorner[\Psi_L\otimes\Psi_R].
\end{eqnarray*}
Next we need a technical lemma taken from~\cite{wi05} Corollary 1.25.

\begin{lem}\label{action}
Let $\alpha$ be a $1$--form. Its metric dual will be also denoted by $\alpha$. Then
\begin{eqnarray*}
\alpha\wedge[\Psi_1\otimes\Psi_2] & = & \phantom{-}\frac{1}{2}\big((-1)^m[\alpha\cdot\Psi_1\otimes\Psi_2]-\widetilde{[\Psi_1\otimes \alpha\cdot\Psi_2]}\big),\\
\alpha\llcorner[\Psi_1\otimes\Psi_2] & = & \frac{1}{2}\big(-(-1)^m[\alpha\cdot\Psi_1\otimes\Psi_2]-\widetilde{[\Psi_1\otimes \alpha\cdot\Psi_2]}\big).
\end{eqnarray*}
\end{lem}
For $H\in\Omega^3(M)$ it follows 
\begin{eqnarray*}
H\wedge[\Psi_1\otimes\Psi_2] & = & \phantom{-}\frac{(-1)^m}{8}\big[H\cdot\Psi_1\otimes\Psi_2-\sum_ke_k\cdot\Psi_1\otimes(e_k\llcorner H)\cdot\Psi_2\big]\\
& & +\frac{1}{8}\big[\Psi_1\otimes H\cdot\Psi_2-\sum_k(e_k\llcorner H)\cdot\Psi_1\otimes e_k\cdot\Psi_2\big]^{\sim},\\
H\llcorner[\Psi_1\otimes\Psi_2] & = &  \phantom{-}\frac{(-1)^m}{8}\big[-H\cdot\Psi_1\otimes\Psi_2+\sum_ke_k\cdot\Psi_1\otimes(e_k\llcorner H)\cdot\Psi_2\big]\\
& & +\frac{1}{8}\big[\Psi_1\otimes H\cdot\Psi_2-\sum_k(e_k\llcorner H)\cdot\Psi_1\otimes e_k\cdot\Psi_2\big]^{\sim}.
\end{eqnarray*}
Since $\nabla$ and $[\cdot\,,\cdot]$ commute (for $\nabla$ is metric), $d=\sum e_k\wedge\nabla_{e_k}$ and $d^*=-\sum e_k\llcorner\nabla_{e_k}$ imply
\begin{eqnarray}
\big[\mc{D}(\Psi_1\otimes\Psi_2)\big] & = & (-1)^m\big(d^*[\Psi_1\otimes\Psi_2]+d[\Psi_1\otimes\Psi_2]\big)\nonumber\\
\big[\widetilde{\mc{D}}(\Psi_1\otimes\Psi_2)\big] & = & \big(d^*[\Psi_1\otimes\Psi_2]-d[\Psi_1\otimes\Psi_2]\big)^{\sim}\label{diracvsd}.
\end{eqnarray}
As a result,
\begin{eqnarray}
\big[\mc{D}\big(\mc{A}(\Psi_L)\otimes\Psi_R\big)\big] & = & (-1)^me^{\phi}(F_0+\widehat{m}i^m\gtilde F_0)\nonumber\\
& & +(-1)^m\big((\alpha+H)\llcorner-(\alpha+H)\wedge\big)[\mc{A}(\Psi_L)\otimes\Psi_R]\nonumber\\
\big[\mc{D}(\Psi_L\otimes\Psi_R)\big] & = & (-1)^me^{\phi}(F_1+\widehat{m}i^m\gtilde F_1)\nonumber\\
& & +(-1)^m\big((\alpha+H)\llcorner-(\alpha+H)\wedge\big)[\Psi_L\otimes\Psi_R]\nonumber
\\
\big[\widetilde{\mc{D}}(\mc{A}(\Psi_L)\otimes\Psi_R)\big] & = & e^{\phi}(F_0-\widehat{m}i^m\gtilde F_0)\nonumber\\ & & -\big((\alpha+H)\llcorner+(\alpha+H)\wedge\big)[\mc{A}(\Psi_L)\otimes\Psi_R]\nonumber\\
\big[\widetilde{\mc{D}}(\Psi_L\otimes\Psi_R)\big] & = & (-1)^me^{\phi}(F_1-\widehat{m}i^m\gtilde F_1)\nonumber\\ & & -(-1)^m\big((\alpha+H)\llcorner+(\alpha+H)\wedge\big)[\Psi_L\otimes\Psi_R].\label{twisteddiracev}
\end{eqnarray}
Using again the lemma to compute the action of $\alpha+H$ on $[\Psi_L\otimes\Psi_R]$ the two first equations of~(\ref{twisteddiracev}) become
\begin{eqnarray}
\D\mc{A}(\Psi_L)\otimes\Psi_R+\sum e_k\cdot\mc{A}(\Psi_L)\otimes\nabla_{e_k}\Psi_R & = & (-1)^me^{\phi}\big[F_0+\widehat{m}i^m\gtilde F_0\big]^{-1}\label{tweq2}\\
& &+\sum e_k\cdot\mc{A}(\Psi_L)\otimes(e_k\llcorner\frac{H}{4})\cdot\Psi_R\nonumber\\
& &-\alpha\cdot\mc{A}(\Psi_L)\otimes\Psi_R-\frac{H}{4} \cdot\mc{A}(\Psi_L)\otimes\Psi_R,\nonumber\\
\D\Psi_L\otimes\Psi_R+\sum e_k\cdot\Psi_L\otimes\nabla_{e_k}\Psi_R & = &(-1)^me^{\phi}\big[F_1+\widehat{m}i^m\gtilde F_1\big]^{-1}-\frac{H}{4}\cdot\Psi_L\otimes\Psi_R\label{tweq1}\nonumber\\
& & +\sum e_k\cdot\Psi_L\otimes(e_k\llcorner\frac{H}{4})\cdot\Psi_R-\alpha\cdot\Psi_L\otimes\Psi_R.\nonumber
\end{eqnarray}
Next we contract the first and the second equation of~(\ref{tweq2}) with $q(e_j\cdot\mc{A}(\Psi_L),\,\cdot)\otimes\mathrm{Id}$ and $q(e_j\cdot\Psi_L,\,\cdot)\otimes\mathrm{Id}$ respectively. Using Lemma~\ref{contlemma} and the identities
\begin{equation}\nonumber
\widehat{\gtilde F_0}=(-1)^{m+1}\gtilde\widehat{F}_0,\quad \widehat{\gtilde F_1}=-\gtilde\widehat{F}_1,
\end{equation}
we obtain
\begin{eqnarray}
\sum_kq\big(e_j\cdot\mc{A}(\Psi_L),e_k\cdot\mc{A}(\Psi_L)\big)\nabla_{e_k}^-\Psi_R & = &-q\big(e_j\cdot\mc{A}(\Psi_L),\D\mc{A}(\Psi_L)+\alpha\cdot\mc{A}(\Psi_L)+\nonumber\\
& &+\frac{1}{4}H\cdot\mc{A}(\Psi_L)\big)\Psi_R-\frac{\widehat{m}}{2^{m-1}}e^{\phi}\widehat{F_0}\cdot e_j\cdot\Psi_L,\label{coneq1}\\
\sum_kq\big(e_j\cdot\Psi_L,e_k\cdot\Psi_L\big)\nabla_{e_k}^-\Psi_R & = &-q\big(e_j\cdot\Psi_L,\D\Psi_L+\alpha\cdot\Psi_L+\frac{1}{4}H\cdot\Psi_L\big)\Psi_R\nonumber\\
& & -\frac{1}{2^{m-1}}e^{\phi}\widehat{F_1}\cdot e_j\cdot\mc{A}(\Psi_L).\label{coneq2}
\end{eqnarray}
Adding~(\ref{coneq1}) and~(\ref{coneq2}) yields
\begin{equation}\nonumber
\nabla_{e_j}^-\Psi_R=-\frac{1}{2^m}e^{\phi}\big(\widehat{m}\widehat{F_0}\cdot e_j\cdot\Psi_L+\widehat{F_1}\cdot e_j\cdot\mc{A}(\Psi_L)\big)-\Re q\big(e_j\cdot\Psi_L,\D\Psi_L+\alpha\cdot\Psi_L+\frac{1}{4}H\cdot\Psi_L\big)\Psi_R,
\end{equation}
for the real part of $q(e_j\cdot\Psi_L,e_k\cdot\Psi_L)$ vanishes unless $j=k$ when it equals $1$. 
We contract once more with $q(\Psi_R,\cdot)$ to see that ${\rm Re}\,q(e_j\cdot\Psi_L,\ldots)=0$ as the remaining terms are purely imaginary or vanish for $F_{0,1}$ is Ramond--Ramond. Hence
\begin{equation}\nonumber
\nabla_X^-\Psi_R=-\frac{1}{2^m}e^{\phi}\big(\widehat{m}\widehat{F_0}\cdot X\cdot\Psi_L+\widehat{F_1}\cdot X\cdot\mc{A}(\Psi_L)\big)
\end{equation}
for all $X\in\mf{X}(M)$. Using the second set of equations in~(\ref{twisteddiracev}) gives similarly
\begin{equation}\nonumber
\nabla_X^+\Psi_L=\frac{1}{2^m}e^{\phi}\big(\widehat{m}\overline{F_0}\cdot X\cdot\Psi_R+m^{\vee}F_1\cdot X\cdot\mc{A}(\Psi_R)\big).
\end{equation}
For the dilatino equations, contracting the second equation in~(\ref{tweq2}) with $\mathrm{Id}\otimes q(\Psi_R,\cdot)$ implies
\begin{equation}\nonumber
\mathrm{D}\Psi_L+q(\Psi_R,\nabla^-_{e_k}\Psi_R)e_k\cdot\Psi_L=-\frac{\widehat{m}}{2^{m-1}}e^{\phi}F_1\cdot\mc{A}(\Psi_R)-\frac{1}{4}H\cdot\Psi_L-\alpha\cdot\Psi_L.
\end{equation}
Contracting the first equation yields (note that $\mathrm{D}\mc{A}=(-1)^{m+1}\mc{A}\mathrm{D}$ etc.)
\begin{equation}\nonumber
\mathrm{D}\Psi_L+\overline{q(\Psi_R,\nabla^-_{e_k}\Psi_R)}e_k\cdot\Psi_L=-\frac{\widehat{m}}{2^{m-1}}e^{\phi}\overline{F_0}\cdot\Psi_R-\frac{1}{4}H\cdot\Psi_L-\alpha\cdot\Psi_L.
\end{equation}
Adding both equations and using that $F_{0,1}$ are Ramond--Ramond fields obtains $\mathrm{D}\Psi_L+H/4\cdot\Psi_L+\alpha\cdot\Psi_L=0$. In the same way one proves the dilatino equation for $\Psi_R$.

\medskip

Let us briefly sketch the converse. By~(\ref{diracvsd}), we have $d[\Psi_L\otimes\Psi_R]=-[\mc{D}(\Psi_L\otimes\Psi_R)+\widetilde{\mc{D}}(\Psi_L\otimes\Psi_R)]/2$. If $F_{0,1}=0$, plugging in the dilatino and the gravitino equation implies in conjunction with Lemma~\ref{action} that $d[\Psi_L\otimes\Psi_R]=(d\phi-H)\wedge[\Psi_L\otimes\Psi_R]$, and similarly for $[\mc{A}(\Psi_L)\otimes\Psi_R]$. In presence of Ramond--Ramond fields, we get additional terms involving $F_0$ and $F_1$. For sake of concreteness, let us assume $m$ to be odd and consider $F_1$. We decompose $F_1=F_{1,+-}\oplus F_{1,-+}$ with $F_{1,ab}\in\Delta_+\otimes\Delta_-$ and write locally 
\begin{eqnarray*}
F_1 & = & F_{1,+-}\oplus F_{1,-+}\\
& = & \sum_{j,l=1}^ma_{jl}\varphi_j\otimes\mc{A}(\psi_l)+\sum_{l=1}^m\sum_rb_{rl}v_{L,r}\otimes\mc{A}(\psi_l)+\\
&  & \sum_{j,l=1}^m\widetilde{a}_{jl}\mc{A}(\varphi_j)\otimes\psi_l+\sum_{j=1}^m\sum_r\widetilde{b}_{jr}\mc{A}(\varphi_j)\otimes v_{R,r}\big),
\end{eqnarray*}
where
\begin{equation}\nonumber
\big(\varphi_j=z_{L,j}\cdot\mc{A}(\Psi_L),\,v_{L,r}\big)\mbox{ and }\big(\psi_j=z_{R,j}\cdot\mc{A}(\Psi_R),\,v_{R,l}\big)
\end{equation}
are orthonormal bases for $\C^m_L\oplus V_L$ and $\C^m_R\oplus V_R$ respectively. By~(\ref{faction}),
\begin{eqnarray*}
\frac{1}{2^m}\sum_{k=1}^{2m}\widehat{F_1}\cdot e_k\cdot\mc{A}(\Psi_L) & = & \sum_{k=1}^{2m}\sum_{j,l=1}^m\widetilde{a}_{jl}q\big(z_{L,j}\cdot\mc{A}(\Psi_L),e_k\cdot\mc{A}(\Psi_L)\big)\psi_l+\\
& &\sum_{k=1}^{2m}\sum_{j=1}^m\sum_r\widetilde{b}_{jr}q\big(z_{L,j}\cdot\mc{A}(\Psi_L),e_k\cdot\mc{A}(\Psi_L)\big)v_{R,r}.
\end{eqnarray*}
Using~(\ref{Jaction}) finally gives
\begin{eqnarray*}
\frac{1}{2^{m+1}}\sum_{k=1}^{2m}e_k\cdot\Psi_L\otimes\widehat{F_1}\cdot e_k\cdot\mc{A}(\Psi_L)&=&\sum_{k,\,j,\,l=1}^{m}\!\widetilde{a}_{jl}\mc{A}(\varphi_k)\otimes q\big(z_{L,j}\cdot\mc{A}(\Psi_L),z_{L,k}\cdot\mc{A}(\Psi_L)\big)\psi_l+\\
& &\sum_{k,\,j,\,r}\widetilde{b}_{jr}\mc{A}(\varphi_k)\otimes q\big(z_{L,j}\cdot\mc{A}(\Psi_L),z_{L,k}\cdot\mc{A}(\Psi_L)\big)v_{R,r}\\
& =& F_{1,-+}.
\end{eqnarray*}
Proceeding similarly with the terms involving $F_0$ concludes the proof.
\end{prf}

\begin{rmk}
As the proof shows, the implication ``$\Rightarrow$'' is valid for any $F_{0,1}$ satisfying~(\ref{rrchar}). Further,  the T--duality mechanism as discussed in~\cite{gmwi06} preserves the condition~(\ref{rrchar}) by equivariance. Since it also preserves $d_{\nu}$--exactness, one can T--dualise an internal space of type IIB space into one of type IIA and vice versa, in agreement with physical expectations.
\end{rmk}

\begin{ex}
We continue the last example of the previous section with $\rho_0=e^{-i\omega}$ and $\rho_1=\psi_++i\psi_-$. Critical points as in (i) are obtained for $F_0=d\omega$ provided $d\psi_{\pm}=0$ and $d\omega\wedge\omega=0$. It remains to see under which conditions $F_0\cdot\Psi=0$ hold, that is, $F_0$ is a Ramond--Ramond field. A routine application of Schur's lemma implies $F_0\in\llbracket\Lambda_0^{2,1}\rrbracket$, following the notation of~\cite{chsa02} -- this automatically implies $d\omega\wedge\omega=0$, cf.~\cite{chsa02}. Still following the same paper, we refer to such an $SU(3)$--structure as being of {\em type $W_3$}. Similarly, $F_1=d\psi_{\pm}$ defines critical points as in (ii) provided $d\omega=0$ and $d\psi_{\mp}=0$. Moreover, $F_1$ is Ramond--Ramond \iff $F_1\in\llbracket\Lambda_0^{1,1}\rrbracket$, that is, $M$ is either of {\em type} $W_2^+$ or $W_2^-$. Explicit examples hereof can be found in~\cite{chsa02}.
\end{ex}

\bigskip

{\em Acknowledgements.} The authors would like to thank Stefano Chiantese, Florian Gmeiner, Dieter L\"ust and Dimitrios Tsimpis for useful discussions. The second author also would like to thank the Erwin Schr\"odinger Institute, Vienna, for hospitality during the preparation of this paper.
%
%
%
%
%

\end{document}